\documentclass[5p,twocolumn,10pt]{elsarticle}
\usepackage{mathrsfs}
\usepackage{mathrsfs}
\usepackage{graphicx}
\usepackage{color}
\usepackage{float}
\usepackage{amsfonts}
\usepackage{graphicx,times,amsmath}
\usepackage{graphicx}
\usepackage{epstopdf}
\usepackage{latexsym}
\usepackage{amsmath}
\usepackage{amsfonts}
\usepackage{amssymb}
\usepackage{url}
\usepackage{subfigure}
\usepackage{times}
\usepackage{color}
\usepackage{amsmath}
\usepackage{amssymb}
\usepackage{lineno}
\newtheorem{theorem}{Theorem}

\newtheorem{corollary}{Corollary}
\newtheorem{definition}{Definition}
\newtheorem{lemma}{Lemma}
\newtheorem{remark}{Remark}
\newtheorem{notation}{Notation}
\newtheorem{assumption}{Assumption}

\allowdisplaybreaks
\journal{Neurocomputing}
\begin{document}

\begin{frontmatter}

\title{Global $\mu$-stability and finite-time control of octonion-valued neural networks with unbounded delays\tnoteref{tt}}\tnotetext[tt]{This work was supported by the National Science Foundation of China under Grant No. 61673298, 61203149; Shanghai Rising-Star Program of China under Grant No. 17QA1404500; Natural Science Foundation of Shanghai under Grant No. 17ZR1445700; the Fundamental Research Funds for the Central Universities of Tongji University.}


\author{Jingzhu Wang}
\author{Xiwei Liu\corref{cor1}}
\address{Department of Computer Science and Technology, Tongji University, and with the Key Laboratory of Embedded System and Service Computing, Ministry of Education, Shanghai 201804, China.}
\cortext[cor1]{Corresponding Author. E-mail address: xwliu@tongji.edu.cn}

\begin{abstract}
Octonion-valued neural networks (OVNNs) are a type of neural networks for which the states and weights are octonions. In this paper, the global $\mu$-stability and finite-time stability problems for octonion-valued neural networks are considered under unbounded and asynchronous time-varying delays. To avoid the non-communicative and non-associative multiplication feature of the octonions, we firstly decompose the OVNNs into eight real-valued neural networks (RVNNs) equivalently. Through the use of generalized norm and the Cauchy convergence principle, we obtain the sufficient criteria which assure the existence, uniqueness of the equilibrium point and global $\mu$-stability of OVNNs. By adding controllers, the  criteria to ensure the finite-time stability for OVNNs are presented by dividing the analysis of finite-time stability process into two phases. Furthermore, we also prove the adaptive finite time stability theory of above networks. At last, the simulation results of specified examples is given to substantiate  the effectiveness and  correctness of the theoretical results.
\end{abstract}


\end{frontmatter}

\section{Introduction}
In the last few years, researches on neural networks (NNs) with values in multidimensional domains have aroused great interest, such as complex-valued neural networks (CVNNs) and quaternion-valued neural networks (QVNNs), for which the states, connection weights and activation functions are complex-valued and  quaternion-valued respectively. A lot of applications of CVNN and QVNN models in different fields have been discovered and because of the requirements of the stability of NNs in most applications, many theoretical analyses about the stability of equilibrium were proposed, see \cite{cvnn1}-\cite{qvnn5}.
	
Recently, it has been natural to extend research on CVNN to OVNN, because the characteristic dimensions of research objects are becoming more and more high-dimensional and octonions are widely used in geometry, physics and signal processing, \cite{app1}-\cite{app4}. As a newly developed multidimensional NNs, OVNN was first proposed by Popa \cite{popa1}. As an extension of CVNN and QVNN, OVNN has octonion-valued states, octonion-valued connection weights, octonion-valued activation functions and octonion-valued output. In addition, octonion algebra has the significant property of normed division algebra, i.e., we can define norm and multiplicative inverse on it. What's more, it can be proved that octonion algebra is the only norm-division algebra that can be defined in the real field except complex algebra and quaternion algebra, see \cite{popa2} .
	
Similar to CVNN, it is necessary to study the stability of OVNN. To the best of our knowledge, there have been few papers about the stability of QVNNs. \cite{popa2} proposed sufficient conditions for the global exponential stability of neutral-type OVNNs with time-varying delays. Utilizing the Cayley-Dickson construction, the OVNNs system was separated into four complex-valued systems. Then based on two Lipschitz condition assumptions, the stability criteria are established in forms of complex-valued linear matrix inequalities(LMIs). In \cite{popa3}, a sufficient criterion was obtained by transforming the octonion-valued differential into eight real-valued	equations and taking advantage of LMIs, which guaranteed the existence, uniqueness and global asymptotic stability of the equilibrium point for OVNNs with delay. In \cite{popa4} and \cite{popa5}, the exponential stability for OVNNs with delay, leakage delay, time-varying delays and distributed delays was discussed.
	
Since there is inevitably delay in practice, papers mentioned above all study the cases with delay or time-varying delay and some stability criteria for delayed OVNN systems were reported. However, it should be noticed that the time delays are assumed to be bounded in most of the papers, which means $\tau_{pq}(t) \leq \tau$, $\tau$ is a fixed constant. In order to deal with the case that there may be systems with unbounded time-delays in real life, \cite{mu1} introduced a new concept of $\mu$-stability, which was a generalization of power-stability discussed in \cite{mu2}. This concept has been applied to CVNNs and QVNNs (see \cite{cvnn5}, \cite{qvnn5}, \cite{mu3}), but not OVNNs yet.

It is worth noting that	above studies are based on the infinite convergence time, which takes a higher cost. In practical engineering applications, systems are preferred to reach stability in finite time, instead of infinite time. In order to realize stability of the system quickly, researchers use finite-time techniques to get a faster convergent speed and finite-time stability problems in dynamical systems has been extensively studied in recent years. What's more, as an important application of finite-time stability, finite-time synchronization/anti-synchronization problem has also attracted more and more attentions. Some control strategies are established and applied in these cases for achieving the finite-time stability/synchronization/anti-synchronization of networks \cite{FTSta1}-\cite{asyn3}. \cite{asyn1} divided the whole anti-synchronization process of master-slave NNs with time delays into two procedures and prove that both are finite. \cite{asyn3} investigated the finite time anti-synchronization of master-slave coupled CVNNs with bounded asynchronous time-varying delays and designed a controller with three control terms. \cite{FTSta3} set up some theories about (adaptive) finite-time stability with unbounded time-varying delays of a general model.

Motivated by the above analysis, in this paper,

(1) The global $\mu$-stability of OVNNs with unbounded time-varying delays will be considered.

Considering that time delays between different neurons are usually various in practice, the time delays we discussed are asynchronous. To avoid the non-communicative and non-associative multiplication feature of the octonion algebra, we firstly decompose the OVNNs into eight real-valued neural networks (RVNNs) equivalently. Through discussion about companion neural networks of studied networks, we utilize generalized norm and the Cauchy convergence principle to get some criteria in order to make sure equilibrium point exist uniquely.

%Besides, the OVNN in this paper is general that obtained criteria can reduce to QVNN in \cite{qvnn5} and CVNN in \cite{cvnn5}.
	
(2) The finite-time	stability of OVNNs with unbounded asynchronous time-varying delays will be investigated.

Once the stability condition is not satisfied, then the external controller should be added. The controller in this paper will be designed as simple as possible, i.e., the number of control terms is reduced to two and is independent from time delays. Using the decomposing technique and generalized norm, sufficient criteria are obtained to ensure the OVNNs realize finite-time stability. We will take the two-phases-method (2PM, \cite{asyn1,FTSta3}) to prove the finite-time stability of OVNNs with time delays. Moreover, we design a simple adaptive rule, such that adaptive finite-time stability can also be realized.
	
The rest of the paper is organized as follows. Section \ref{model} introduces some basic properties of octonion algebra and gives the definition of QVNNs with unbounded and asynchronous time-varying delays, together with some necessary definitions, assumptions and lemmas. In section \ref{theory}, we derive some sufficient criteria based on generalized norm to ensure the global $\mu$-stability and (adaptive) finite-time stability of considered OVNNs. In section \ref{simu}, numerical examples are given to demonstrate the effectiveness of derived results. Finally, some conclusions are drawn in Section \ref{con}.
	
\begin{notation}\label{notation1}
$\mathbb{R}, \mathbb{O}$ represent the sets of real numbers and octonion numbers. $\mathbb {R}$$^{m\times n}$, $\mathbb {O}$$^{m\times n}$ represent any ${m \times n}$-dimensional real-valued and  octonion-valued matrices.
\end{notation}
	
\begin{notation}\label{notation2}
$\{a\}^+=\max\{0,a\}$ when $a \in \mathbb{R}$. For any vector $a=(a_1,\cdots,a_n)$, then $\{a\}^+=(\{a_1\}^+,\cdots, \{a_n\}^+)$, $|a|=(|a_1|,\cdots,|a_n|)$, $a^m=(a_1^m,\cdots,a_n^m)$.
	\end{notation}
	
\section{Problem model and preliminaries}\label{model}
At first, the definition and some properties of octonion algebra are introduced.
	
The octonions which we denote $\mathbb {O}$ are an 8-dimensional algebra with basis	$\{e_0, e_1, \cdots, e_7\}$. The algebra of octonions consists of elements such as $x=\sum_{\ell=0}^7 x^{\ell}e_{\ell}, \ell=0,1,\cdots,7$, where $x^{\ell} \in \mathbb {R}$ is a real coefficient and $e_{\ell}$ is the fundamental octonion unit  satisfying the following multiplication rule table:
\begin{center}
	\begin{tabular}{| c | c | c | c | c | c | c | c | c |}
		\hline
		$\times $ & $e_0$ & $e_1$ & $e_2$ & $e_3$ & $e_4$ & $e_5$ & $e_6$ & $e_7$ \\
		\hline
		$e_0$ & $e_0$ & $e_1$ & $e_2$ & $e_3$ & $e_4$ & $e_5$ & $e_6$ & $e_7$ \\
		\hline
		$e_1$ &$e_1$ &$-e_0$ &$e_3$ &$-e_2$ &$e_5$ &$-e_4$ &$-e_7$ &$e_6$ \\
		\hline
		$e_2$ &$e_2$ &$-e_3$ &$-e_0$ &$e_1$ &$e_6$ &$e_7$ &$-e_4$ &$-e_5$ \\
		\hline
		$e_3$ &$e_3$ &$e_2$ &$-e_1$ &$-e_0$ &$e_7$ &$-e_6$ &$e_5$ &$-e_4$ \\
		\hline
		$e_4$ &$e_4$ &$-e_5$ &$-e_6$ &$-e_7$ &$-e_0$ &$e_1$ &$e_2$ &$e_3$ \\
		\hline
		$e_5$ &$e_5$ &$e_4$ &$-e_7$ &$e_6$ &$-e_1$ &$-e_0$ &$-e_3$ &$e_2$ \\
		\hline
		$e_6$ &$e_6$ &$e_7$ &$e_4$ &$-e_5$ &$-e_2$ &$e_3$ &$-e_0$ &$-e_1$ \\
		\hline
		$e_7$ &$e_7$ &$-e_6$ &$e_5$ &$e_4$ &$-e_3$ &$-e_2$ &$e_1$ &$-e_0$\\
		\hline
	\end{tabular}
\end{center}
	
The \emph{addition of octonions} is defined by $a+b=\sum_{\ell=0}^7(a^\ell+b^\ell)e_{\ell}$, where $a=\sum_{\ell=0}^7a^\ell e_{\ell}$, $b=\sum_{\ell=0}^7b^\ell e_{\ell}$. \emph{Scalar multiplication} is defined as $ma=\sum_{\ell=0}^7(ma^\ell)e_\ell$.

An \emph{octonion function} $f(t) \in \mathbb {O}$ is described by $f(t)=\sum_{\ell=0}^7f^\ell(t) e_{\ell}$, where $f^\ell(t)$ are real-valued functions.
	
Now, the properties of the product between any	octonion-valued numbers $a$, $b$ are introduced. Firstly we define a multiplication matrix	
\begin{align}
&M=\sum_{\ell=0}^7M^{\ell}e_{\ell}\label{equation1}\\
=& \begin{pmatrix}
	e_0 &e_1  &e_2 &e_3 &e_4 &e_5  &e_6  &e_7 \\
	e_1 &-e_0 &e_3 &-e_2 &e_5 &-e_4 &-e_7 &e_6 \\
	e_2 &-e_3 &-e_0 &e_1 &e_6 &e_7 &-e_4 &-e_5 \\
	e_3 &e_2 &-e_1 &-e_0 &e_7 &-e_6 &e_5 &-e_4 \\
	e_4 &-e_5 &-e_6 &-e_7 &-e_0 &e_1 &e_2 &e_3 \\
	e_5 &e_4 &-e_7 &e_6 &-e_1 &-e_0 &-e_3 &e_2 \\
	e_6 &e_7 &e_4 &-e_5 &-e_2 &e_3 &-e_0 &-e_1 \\
	e_7 &-e_6 &e_5 &e_4 &-e_3 &-e_2 &e_1 &-e_0
	\end{pmatrix} 	\nonumber
	\end{align}
where matrices $M^0=\mathrm{diag}(1,-1,-1,-1,-1,-1,-1,-1)$;
\begin{align*}M^1=
	\begin{pmatrix}
	0 &1 &0 &0 &0 &0 &0 &0 \\
	1 &0 &0 &0 &0 &0 &0 &0 \\
	0 &0 &0 &1 &0 &0 &0 &0 \\
	0 &0 &-1 &0 &0 &0 &0 &0 \\
	0 &0 &0 &0 &0 &1 &0 &0 \\
	0 &0 &0 &0 &-1 &0 &0 &0 \\
	0 &0 &0 &0 &0 &0 &0 &-1 \\
	0 &0 &0 &0 &0 &0 &1 &0 	
	\end{pmatrix}
\end{align*}
Obviously, there are many zeros in the above matrix, and due to the limit of page, we represent the above matrix as:
$M^1(1,2)=M^1(2,1)=M^1(3,4)=M^1(5,6)=M^1(8,7)=1$, $M^1(4,3)=M^1(6,5)=M^1(7,8)=-1$, and other not mentioned elements are zero.

Similarly, $M^2(1,3)=M^2(3,1)=M^2(4,2)=M^2(5,7)=M^2(6,8)=1$, $M^2(2,4)=M^2(7,5)=M^2(8,6)=-1$; $M^3(1,4)=M^3(2,3)=M^3(4,1)=M^3(5,8)=M^3(7,6)=1$, $M^3(3,2)=M^3(6,7)=M^3(8,5)=-1$; $M^4(1,5)=M^4(5,1)=M^4(6,2)=M^4(7,3)=M^4(8,4)=1$, $M^4(2,6)=M^4(3,7)=M^4(4,8)=-1$; $M^5(1,6)=M^5(2,5)=M^5(6,1)=M^5(8,3)=1$, $M^5(3,8)=M^5(4,7)=M^5(5,2)=M^5(7,4)=-1$; $M^6(1,7)=M^6(2,8)=M^6(3,5)=M^6(6,4)=M^6(7,1)=1$, $M^6(4,6)=M^6(5,3)=M^6(8,2)=-1$; $M^7(1,8)=M^7(3,6)=M^7(4,5)=M^7(7,2)=M^7(8,1)=1$, $M^7(2,7)=M^7(5,4)=M^7(6,3)=-1$.

	%delete space before lemma.dont know why
	%\par\setlength\parindent{0em}
\begin{lemma}\label{lemma1}
For any two numbers $a$, $b\in\mathbb {O}$, define
\begin{align}\label{equation2}
		\tilde{a}=(a^0,a^1,\cdots,a^7)^T, ~~ \tilde{b}=(b^0,b^1,\cdots,b^7)^T
\end{align}
then $a b =\tilde{a}^TM\tilde{b}=\sum_{\ell=0}^7e^{\ell}\tilde{a}^TM^{\ell}\tilde{b}$,
or, $\tilde{ab}=(\tilde{a}^TM^0\tilde{b}, \tilde{a}^TM^1\tilde{b}, \cdots, \tilde{a}^TM^7\tilde{b})^T$.

	\end{lemma}
	
Its proof is based on the definition of $\tilde{a}$ and $\tilde{b}$, and making use of the multiplication rules of octonion units,
\begin{align*}
a\cdot b =&\tilde{a}^T(e_0, e_1,\cdots,e_7)^T\cdot (e_0,e_1,\cdots, e_7) \tilde{b}=\tilde{a}^TM\tilde{b}
\end{align*}	
	
Consider the following OVNN with unbounded time delay:
\begin{align}
		\dot{w}_p(t)=-d_pw_p(t)&+\sum_{q=1}^na_{pq}f_q(w_q(t)) \label{equation5}\\
		&+\sum_{q=1}^nb_{pq}g_q(w_q(t-\tau_{pq}(t)))+I_p,\nonumber
\end{align}
where $w_p(t)\in \mathbb{O}, p=1,\cdots,n$ is the state of $p$-th neuron at time $t$. $D=\mathrm{diag}(d_1,\cdots,d_n) \in \mathbb {R}^{n \times n}$ with $d_p >0$, is the feedback self-connection weight matrix. $a_{pq}\in \mathbb{O}$ and $b_{pq} \in \mathbb{O}$ represent the connection weight between $p$-th neuron and $q$-th neuron without and with time delays, $q=1,\cdots,n$. Functions $f_q(\cdot):\mathbb{O} \rightarrow \mathbb{O}$ and $g_q(\cdot):\mathbb{O} \rightarrow \mathbb{O}$ denote activation functions without and with time delays. $\tau_{pq}(t)$ is the unbounded, asynchronous and time-varying transmission delay between $p$-th and $q$-th neurons. $I_p \in \mathbb{O}$ represents the $p$-th external input.
	
\begin{assumption}\label{ass1}
Suppose octonion function can be decomposed as: $f_q(w_q)=\sum_{\ell=0}^7e_{\ell}f_q^{\ell}(w_q^0,w_q^1,\cdots,w_q^7)$,
where $w_q =\sum_{\ell=0}^7e_{\ell}q^{\ell}$, $f_q^\ell:\mathbb{R}^8 \rightarrow \mathbb{R}, q=1,\cdots,n$. Similar assumptions hold for functions $g_q(w_q)$. Assume the partial derivatives of $f_q(w_q)$ ($g_q(w_q)$) with respect to (w.r.t.) $w_q^0,\cdots,w_q^7$ exist and they are continuous and bounded, i.e., there exist positive constant numbers $\lambda^{\ell_1\ell_2}$ ($\delta^{\ell_1\ell_2}$), $\ell_1,\ell_2=0,1,\cdots,7$, such that
\begin{align*}
0\leq \dfrac{\partial{f_q^{\ell_1}(\cdot)}}{\partial{w_q^{\ell_2}}} \leq \lambda_q^{\ell_1\ell_2},	~~~~
0\leq \dfrac{\partial{g_q^{\ell_1}}(\cdot)}{\partial{w_q^{\ell_2}}} \leq \delta_q^{\ell_1\ell_2}.
			\end{align*}
	\end{assumption}
	
	\begin{assumption}\label{ass2}
For the unbounded time delay, assume we can find a continuous function $\tau(t)$ such that $\tau_{pq}(t) \leq \tau(t)$, where $t-\tau(t)\rightarrow +\infty$ and $\tau(t) \rightarrow +\infty$ as $t \rightarrow +\infty$.
	\end{assumption}
	
\begin{assumption}\label{ass3}
For above $\tau(t)$, assume $\mu(t)$ is a nondecreasing and continuous function such that $\lim_{t \rightarrow \infty}\mu(t)=+\infty$,
\begin{align}\label{ab-equation}
\lim\limits_{t\rightarrow\infty}\frac{\dot{\mu}(t)}{\mu(t)}=\alpha\ge 0, ~~
\lim\limits_{t\rightarrow\infty}\frac{\mu(t)}{\mu(t-\tau(t))}=1+\beta\ge 1.
\end{align}
\end{assumption}
	
For convenience, we denote $f_q^\ell=f_q^\ell(w_q(t))$, $g_{q\tau_{pq}}^\ell=g_q^\ell(w_q(t-\tau_{pq}(t)))$, $\ell=0,1,\cdots,7, p,q=1,\cdots,n$. According to Assumption \ref{ass1} and Lemma \ref{lemma1}, OVNN (\ref{equation5}) can be separated into eight RVNNs as
\begin{align*}
%\label{equation6}	
\dot{w}_p^{\ell}(t)=-d_pw_p^{\ell}+\sum\limits_{q=1}^n\tilde{a}_{pq}^TM^{\ell}\tilde{f}_q+\sum\limits_{q=1}^n\tilde{b}_{pq}^TM^{\ell}\tilde{g}_{q\tau_{pq}}+I_p^{\ell}
\end{align*}
	
\begin{lemma}\label{lemma2}
Using the partial derivative of $f_q(w_q(t))$, we define a $\mathbb{R}^{8\times 8}$ matrix
		\begin{align*}
		\overline{M}(f_q(w_q(t)))=&\bigg(\frac{\partial f_q^{\ell_1}}{\partial w_q^{\ell_2}}\bigg)
		=&\sum\limits_{c=0}^7\overline{M}(f_q(w_q(t)),c),
		\end{align*}
where $\ell_1, \ell_2=0,1,\cdots,7$, $\overline{M}(f_q(w_q(t)),c)$, $c=0,1,\cdots,7$ means a $8\times 8$ matrix with the $c$-th column in $\overline{M}(f_q(w_q(t)))$, and elements in the other columns are all zeros.
		
Assume $a \in\mathbb{O}$ is a constant, we have
		$$		\frac{d}{dt}\tilde{a}^TM^{c}\widetilde{f_q}(w_q(t))=\tilde{a}^TM^{c}\overline{M}(f_q(w_q(t)))\widetilde{\dot{w}}_q(t)
		$$
	\end{lemma}
	
To simplify expression in the following analysis, by using the upper bounds of the partial derivatives in Assumption \ref{ass1}, we can define the following $\mathbb{R}^{8\times 8}$ matrices
	\begin{align*}
&\overline{M}(q,\lambda)=(\lambda_q^{\ell_1\ell_2}),
~~\overline{M}(q,\delta)=(\delta_q^{\ell_1\ell_2})
	\end{align*}
where $\ell_1, \ell_2=0,1,\cdots,7$,  $\overline{M}_c(q,\lambda)$ (or $ \overline{M}_c(q,\delta)$) represents a vector whose elements are  totally the same with the $c$-th column in  $\overline{M}(q,\lambda)$ (or $ \overline{M}(q,\delta)$).
	
\begin{definition} (\cite{mu2})\label{def1}
For vector $z(t)=(z_1(t),\cdots,z_m(t))^T\in\mathbb{R}^m$, its $\{\xi,\infty\}$-norm is denoted as:
\begin{align*}
		\|z(t)\|_{\{\xi,\infty\}}=\max_j|\xi_j^{-1}z_j(t)|,
\end{align*}
where $\xi=(\xi_1,\cdots,\xi_m)^T>0$. Particularly, when $\xi=(1,\cdots,1)^T$, $\lVert z(t)\rVert_{\{\xi,\infty\}}=\lVert z(t)\rVert_{\infty}$.
	\end{definition}
	
	\begin{definition}(\cite{mu1})\label{def2}
Assume $\mu(t)$ is a nonnegative and continuous function defined in Assumption \ref{ass3}, and $z^*$ is an equilibrium point of NN (\ref{equation5}). If there exists a scalar $W>0$ such that
		$$
		\lVert w(t)-z^*\rVert_{\{\xi,\infty\}}\leq \frac{W}{\mu(t)}, ~~~t\ge 0,
		$$
then NN (\ref{equation5}) is regarded to have global $\mu$-stability.
	\end{definition}
	
	\begin{definition} (\cite{FTSta1})\label{def3}
The equilibrium $z^*$ is said to be a finite time stable equilibrium if the finite time convergence condition and Lyapunov stability condition hold globally.
	\end{definition}
	
\section{Main results}\label{theory}
	
\subsection{Existence and uniqueness of the equilibrium point}
In this subsection, we give a theorem to guarantee that the equilibrium point of OVNN (\ref{equation5}) exists and is unique.

At first, we introduce the following notations:
\begin{align}
\Lambda=&(\xi_1,\cdots,\xi_n,\phi_1,\cdots,\phi_n,\eta_1,\cdots,\eta_n,\zeta_1,\cdots,\zeta_n,\nonumber\\
&~\sigma_1,\cdots,\sigma_n,\psi_1,\cdots, \psi_n,\nu_1,\cdots,\nu_n,\rho_1,\cdots,\rho_n)^T\label{lam}\\
\Lambda_p=&(\xi_p,\phi_p,\eta_p,\zeta_p,\sigma_p,\psi_p,\nu_p,\rho_p)^T,~~p=1,\cdots,n\label{lamp}
\end{align}
and $\Lambda_p^{[\ell]}$ is a vector whose elements are the same with $\Lambda_p$ except the $\ell$-th element is zero, and $\Lambda_p[\ell]$ means the $\ell$-th element of $\Lambda_p$. Now, we have the following result.
\begin{theorem}\label{theorem1}
Suppose Assumption \ref{ass1} holds, OVNN (\ref{equation5}) has a unique equilibrium point if there exists a vector $\Lambda>0$ in (\ref{lam}), such that, for any $p=1,\cdots,n$ and $\ell=0,1,\cdots,7$,
\begin{align} \label{t0}
		T_{\ell}(p)
		=&\Lambda_p[\ell]\Big[-d_p+\{\tilde a_{pp}^T  M^{\ell}\}^+\overline{M}_{\ell}(p,\lambda) \nonumber\\
		&+\{\tilde b_{pp}^T M^{\ell}\}^+\overline{M}_{\ell}(p,\delta)\Big] \nonumber\\
		&+\Big[|\tilde a_{pp}^T M^{\ell}| \overline{M}(p,\lambda)+|\tilde b_{pp}^T M^{\ell}| \overline{M}(p,\delta)\Big]\Lambda_p^{[\ell]} \nonumber\\
		&+\sum_{q\ne p}\Big[|\tilde a_{pq}^T M^{\ell}| \overline{M}(q,\lambda)+|\tilde b_{pq}^T M^{\ell}|  \overline{M}(q,\delta)\Big]\Lambda_q \nonumber\\
		&<0.
\end{align}

	\end{theorem}	
Its proof can be found in Appendix A.

Taking no account of the sign, the following result can be naturally obtained by using the similar method in Theorem \ref{theorem1}. Hence, Theorem \ref{theorem1} has a wilder range of arguments.
\begin{corollary}
Suppose Assumption \ref{ass1} satisfies, OVNN (\ref{equation5}) has a unique equilibrium point if there exists a vector $\Lambda>0$, such that, for any $p=1,\cdots,n$ and $\ell =0,1,\cdots,7$
		\begin{align*}
T_{\ell}(p)=&-d_p\cdot  \Lambda_p[\ell]\\
&+\sum\limits_{q=1}^n\Big[|\tilde{a}_{pq}^TM^{\ell}|\overline{M}(q,\lambda)+|\tilde{b}_{pq}^TM^{\ell}|\overline{M}(q,\delta)\Big]\Lambda_q<0.\nonumber
		\end{align*}
\end{corollary}

\subsection{$\mu$-stability of OVNN}
Now, we will consider the stability of the unique equilibrium $Z^*=(z_1^*,z_2^*,\cdots,z_n^*)\in \mathbb{O}^n$ for delayed OVNN (\ref{equation5}).

Let $\overline{w}_p(t)=w_p(t)-z_p^*=\sum_{\ell=0}^7e_{\ell}\overline{w}_p^{\ell}(t)$, where $w_p(t)$ is the solution of (\ref{equation5}) and $\overline{w}_p^{\ell}(t)=w_p^{\ell}(t)-{z_p^*}^{\ell}, p=1,\cdots,n$.
		
Therefore,
\begin{align*}
%\label{equation*}			
\dot{\overline{w}}_p(t)=-d_p\overline{w}_p(t)+\sum_{q=1}^n\Big[a_{pq}\overline{f}_q(\overline{w}_q)+b_{pq}\overline{g}_{q\tau_{pq}}(\overline{w}_{q\tau_{pq}})\Big],
\end{align*}
where $\overline{f}_q(\overline{w}_q)$=$f_q(\overline{w}_q(t)+z_q^*)-f_q(z_q^*)$ and $\overline{g}_{q\tau_{pq}}(\overline{w}_{q\tau_{pq}})$=\\$g_q(\overline{w}_{q\tau_{pq}}+z_q^*)-g_q(z_q^*)$. For the sake of simplicity, we denote $f_q(z_q^*)$ and $g_q(z_q^*)$ as $f_q^*$ and $g_q^*$ respectively. Then the above equation can be decomposed as follows :
\begin{align*}
		\dot{\overline w}_p^{\ell}(t)=&-d_p\overline{w}_p^{\ell}(t)\\		&+\sum\limits_{q=1}^n\Big[\tilde{a}_{pq}^TM^{\ell}(\tilde{f}_q-\tilde{f}_q^*)+\tilde{b}_{pq}^TM^{\ell}(\tilde{g}_{q\tau_{pq}}-\tilde{g}_q^*)\Big]
		\end{align*}

Denote
		\begin{align} Y(t)=&(\overline{w}_1^0,\cdots,\overline{w}_n^0,\overline{w}_1^1,\cdots,\overline{w}_n^2,\overline{w}_1^3,\cdots,\overline{w}_n^3,\overline{w}_1^4,\cdots,\overline{w}_n^4,\nonumber\\
&~\overline{w}_1^5,\cdots,\overline{w}_n^5,\overline{w}_1^6,\cdots,\overline{w}_n^6,\overline{w}_1^7,\cdots,\overline{w}_n^7,\overline{w}_1^8,\cdots,\overline{w}_n^8)^T\label{yt}
		\end{align}
with
\begin{align}\label{normyt}
\lVert Y(t)\rVert_{\{\Lambda,\infty\}}=\max_{\ell}\Big\{\max_p\{|\Lambda_p[\ell]^{-1}\overline{w}_p^{\ell}(t)|\}\Big\}
\end{align}
Then the stability of the equilibrium for OVNN (\ref{equation5}) is equivalent to consider the stability of $\lVert Y(t)\rVert_{\{\Lambda,\infty\}}$.

\begin{theorem}\label{theorem2}
Suppose Assumptions \ref{ass1}, \ref{ass2} and \ref{ass3} work. OVNN (\ref{equation5}) can reach global $\mu$-stability if there exists a vector $\Lambda>0$ in (\ref{lam}), such that the following inequalities hold:
\begin{align} \label{tt0}
&~~~~~~\overline{T}_{\ell}(p, \alpha, \beta) \nonumber\\
&=\Lambda_p[\ell][-d_p+\alpha+\{\tilde a_{pp}^T M^{\ell}\}^+\cdot\overline{M}_{\ell}(p,\lambda)]\nonumber\\
&+|\tilde a_{pp}^T M^{\ell}| \overline{M}(p,\lambda)\Lambda_p^{[\ell]} \nonumber\\
&+\bigg[\sum_{q\ne p}|\tilde a_{pq}^TM^{\ell}| \overline{M}(q,\lambda)+(1+\beta)\sum_{q=1}^{n}|\tilde b_{pq}^T M^{\ell}|\overline{M}(q,\delta) \bigg]\Lambda_q\nonumber\\
&<0,
\end{align}
where $\alpha$ and $\beta$ are defined in (\ref{ab-equation}).

\end{theorem}
Its proof can be found in Appendix B.

Next, we present two corollaries for specific time delays.
\begin{corollary}\label{corollary-exponential stability}
If time delay $\tau_{pq}(t) \leq \tau$, we can choose $\mu(t)=e^{\alpha t}$, where $\alpha$ is a small positive constant, so $\beta=e^{\alpha \tau}-1\geq 0$ in (\ref{ab-equation}).
If there exists a positive vector $\Lambda>0$ in (\ref{lam}) to make inequalities in (\ref{tt0}):
$\overline{T}_{\ell}(p,\alpha,e^{\alpha\tau}-1)<0$ hold, then the equilibrium in OVNN (\ref{equation5}) can reach global exponential stability.
\end{corollary}

\begin{corollary}\label{corollary-power stability}
If time delay $\tau_{pq}(t) \leq \omega t$, $0<\omega<1$, we can choose $\mu(t)=t^{\gamma}$, where $\gamma$ is a small positive constant, so $\alpha=0$ and $\beta=(1-\omega)^{-\gamma}-1$. If there exists a positive vector
$\Lambda>0$ in (\ref{lam}) to make inequalities in (\ref{tt0}):
$\overline{T}_{\ell}(p,0,(1-\omega)^{-\gamma}-1)<0$ hold, then OVNN (\ref{equation5}) can reach global power-stability.
\end{corollary}

\begin{remark}
	More forms of time delay $\tau(t)$ and the related stability rate $\mu(t)$ can be found in \cite{mu1}. 	
\end{remark}
	
\subsection{Finite-time control for OVNN}
Next, we consider the opposite case, i.e., the OVNN (\ref{equation5}) now is not stable. In this case, from Theorem \ref{theorem2}, there must exist at less one $\ell$, such that $\overline{T}_{\ell}(p, \alpha, \beta)>0$. We will design a controller to make OVNN reach any target point
$\hat{Z}^*=(\hat{z}_1^*,\hat{z}_2^*,\cdots,\hat{z}_n^*)\in \mathbb{O}^n$ in finite time. Denote $\hat{w}_p(t)=w_p(t)-\hat{z}_p^*=\sum_{\ell=0}^7e_{\ell}\hat{w}_p^{\ell}(t)=\sum_{\ell=0}^7e_{\ell}[w_p^{\ell}(t)-(\hat{z}_p^*)^{\ell}], p=1,\cdots,n$.
\begin{align}
\dot{\hat w}_p^{\ell}(t)=&-d_p\hat{w}_p^{\ell}(t)+\sum\limits_{q=1}^n\Big[\tilde{a}_{pq}^TM^{\ell}(\tilde{f}_q-\tilde{f}_q(\hat{z}_q^*))\nonumber\\
&+\tilde{b}_{pq}^TM^{\ell}(\tilde{g}_{q\tau_{pq}}-\tilde{g}_q(\hat{z}_q^*))\Big]+U_p^{\ell}(t)+\hat{I}_p^{\ell},\label{ovnn with u}
\end{align}
where
\begin{align}\label{i}
%\hat{I}_p^{\ell}=I_p^{\ell}-d_pz_p^{\ell}+\sum\limits_{q=1}^n\big[\tilde{a}_{pq}^TM^{\ell}\tilde{f}_q(\hat{z}_q)+\tilde{b}_{pq}^TM^{\ell}\tilde{g}_q(\hat{z}_q)\big],
\hat{I}_p^{\ell}=I_p^{\ell}-d_p(\hat{z}_p^*)^{\ell}+\sum\limits_{q=1}^n\big[\tilde{a}_{pq}^TM^{\ell}\tilde{f}_q(\hat{z}_q^*)+\tilde{b}_{pq}^TM^{\ell}\tilde{g}_q(\hat{z}_q^*)\big],
\end{align}
and $U_p^{\ell}(t)$ is control scheme only depending on the system state at present time, which has the form as
\begin{align}\label{controller} U_p^{\ell}(t)=-\mathrm{sign}(\hat{w}_p^{\ell}(t))\big(\kappa_p^{\ell}|\hat{w}_p^{\ell}(t)|+\hat{\kappa}_p^{\ell}\big)
\end{align}
where $\ell=0,1,\cdots,7$, $p=1,\cdots,n$ and $\kappa_p^{\ell},\hat{\kappa}_p^{\ell}$ are positive.

Denote
\begin{align} \hat{Y}(t)=&(\hat{w}_1^0,\cdots,\hat{w}_n^0,\hat{w}_1^1,\cdots,\hat{w}_n^2,\hat{w}_1^3,\cdots,\hat{w}_n^3,\hat{w}_1^4,\cdots,\hat{w}_n^4,\nonumber\\
&~\hat{w}_1^5,\cdots,\hat{w}_n^5,\hat{w}_1^6,\cdots,\hat{w}_n^6,\hat{w}_1^7,\cdots,\hat{w}_n^7,\hat{w}_1^8,\cdots,\hat{w}_n^8)^T\label{hyt}
		\end{align}
with
\begin{align}\label{normhyt}
\lVert \hat{Y}(t)\rVert_{\{\Lambda,\infty\}}=\max_{\ell}\Big\{\max_p\{|\Lambda_p[\ell]^{-1}\hat{w}_p^{\ell}(t)|\}\Big\}
\end{align}
Then the finite time control stability of OVNN (\ref{equation5}) to $\hat{Z}^*$ is equivalent to consider the stability of $\lVert \hat{Y}(t)\rVert_{\{\Lambda,\infty\}}$.

\begin{theorem}\label{theorem3}
Suppose Assumption \ref{ass1}, \ref{ass2} and \ref{ass3} hold, if there exists a vector $\Lambda>0$ in (\ref{lam}), such that for any $p$ and $\ell$,
\begin{align}\label{condition1}
		\kappa_{p}^{\ell}>\Lambda_{p}[\ell]^{-1}\overline{T}_{\ell}(p, \alpha, \beta)
\end{align}
and
		\begin{align}\label{condition2}
		\hat{\kappa}_p^{\ell} > \sum_{q=1}^n \big|\tilde b_{pq}^T  M^{\ell}\big|  \overline{M}(q,\delta)\Lambda_q+|\hat{I}_p^{\ell}|
		\end{align}
then the system (\ref{ovnn with u}) can acquire finite-time stability under the controllers (\ref{controller}).
		
	\end{theorem}

The method to deal with delay for finite time stability used here is 2PM \cite{FTSta3}, concrete proof can be found in Appendix C.
	\begin{remark}
		When $\hat{Z}^*$ is an equilibrium for network (\ref{equation5}), then $\hat{I}_p^{\ell}$ defined in (\ref{i}) will be zero.
	\end{remark}

	\begin{remark}
		When $\Lambda$ is determined, the larger the value of control parameters $\kappa_p^{\ell}$ and $\hat{\kappa}_p^{\ell}$ are, the faster the speed of stabilization is. So once the parameters of OVNNs are unknown, in order to make sure the finite time stability, we will apply the adaptive technique \cite{FTSta3} on control parameters $\kappa_p^{\ell}$ and $\hat{\kappa}_p^{\ell}$.
	\end{remark}

	\begin{theorem}\label{theorem4}
For system (\ref{ovnn with u}) with the following controller:
		\begin{align}\label{new-controller}
			U_p^{\ell}(t)=-\mathrm{sign}(\hat{w}_p^{\ell}(t))\big(\kappa(t)|\hat{w}_p^{\ell}(t)|+\hat{\kappa}(t)\big)
		\end{align}
		with the adaptive rule
		\begin{align}\label{rules of chi}
			\dot{\hat{\kappa}}(t)=
			\left \{
			 \begin{array}{ll}
				0, &\sup_{t-\tau(t)\leq s \leq t}\lVert \hat{Y}(s) \rVert_{\{\Lambda,\infty\}} >1\\
				c_1,&\sup_{t-\tau(t)\leq s \leq t}\lVert \hat{Y}(s) \rVert_{\{\Lambda,\infty\}} \leq 1 \\
				0,  &\sup_{t-\tau(t)\leq s \leq t}\lVert \hat{Y}(s) \rVert_{\{\Lambda,\infty\}} = 0
			\end{array}
			\right .
		\end{align}
		and
		\begin{align}
&~~~\dot{\kappa}(t)=\label{rules of kappa}\\
&
			\left \{
			\begin{array}{ll}
			c_2\mu(t)\lVert \hat{Y}(t) \rVert_{\{\Lambda,\infty\}}, &\sup_{t-\tau(t)\leq s \leq t}\lVert \hat{Y}(s) \rVert_{\{\Lambda,\infty\}} >1\\
			c_3\lVert \hat{Y}(t) \rVert_{\{\Lambda,\infty\}}, &\sup_{t-\tau(t)\leq s \leq t}\lVert \hat{Y}(s) \rVert_{\{\Lambda,\infty\}} \leq 1 \\
			\end{array}
			\right .\nonumber
		\end{align}
where $\ell=0,1,\cdots,7$, $p=1,\cdots,n$, $\Lambda>0$ is defined in (\ref{lam}), $\hat{Y}(t)$ and $\lVert \hat{Y}(t) \rVert_{\{\Lambda,\infty\}}$ are defined in (\ref{hyt}) and (\ref{normhyt}); parameters $c_1,c_2,c_3$ are any positive scalars and the initial values of $\kappa(0)$ and $\hat{\kappa}(0)$ can be any values, here we choose them as $0$. Then the adaptive finite-time stability of system (\ref{ovnn with u}) can be realized.
	\end{theorem}

We still use 2PM to solve this problem. Its concrete proof can be found in Appendix D.
	
	\begin{corollary}
		Let the $e_4-e_7$ parts of OVNN (\ref{equation5}) and (\ref{ovnn with u}) fixed to be zero and replace the $e_0-e_3$ units with $1,i,j,k$, then the octonion-valued states of OVNN (\ref{equation5}) and (\ref{ovnn with u}) can be transformed to quaternion-valued states, i.e., Theorem \ref{theorem1}-\ref{theorem4} can be applied to QVNNs. On the basis of QVNN, only units R and I are retained, then Theorem \ref{theorem1}-\ref{theorem4} can make available for CVNNs, even RVNNs. Therefore,	results in this paper have a wide scope.
	\end{corollary}
	
	\section{Numerical simulations}\label{simu}

\subsection{$\mu$-stability of OVNN with delay}
Consider a two-neuron OVNN described as
\begin{align}\label{numerical-ovnn}
\left\{
\begin{array}{ll}
\dot{w}_1(t)=&-d_1w_1(t)+a_{11}f_1(w_1(t))\\
&+a_{12}f_2(w_2(t))+b_{11}g_1(w_1(t-\tau_{11}(t)))\\
&+b_{12}g_2(w_2(t-\tau_{12}(t)))+I_1,\\
\dot{w}_2(t)=&-d_2w_2(t)+a_{21}f_1(w_1(t))\\
&+a_{22}f_2(w_2(t))+b_{21}g_1(w_1(t-\tau_{21}(t)))\\
&+b_{22}g_2(w_2(t-\tau_{22}(t)))+I_2,
\end{array}
\right.
\end{align}
where $w_p=\sum_{\ell=0}^7e_{\ell}w_p^{\ell}\in \mathbb{O},p=1,2$, $d_1=d_2=30$ and
\begin{align*}
&\tilde{a}_{11}=(-0.2,-0.3,-0.2,-0.3,-0.2,-0.3,-0.2,-0.3)^T;\\
&\tilde{a}_{12}=(0.3,-0.1,0.3,-0.1,0.3,-0.1,0.3,-0.1)^T;\\
&\tilde{a}_{21}=(0.4,-0.2,0.4,-0.2,0.4,-0.2,0.4,-0.2)^T;\\
&\tilde{a}_{22}=(-0.1,0.2,-0.1,0.2,-0.1,0.2,-0.1,0.2)^T;\\
&\tilde{b}_{11}=(-0.11,0.12,-0.11,0.12,-0.11,0.12,-0.11,0.12)^T;\\
&\tilde{b}_{12}=(0.12,0.11,0.12,0.11,0.12,0.11,0.12,0.11)^T;\\
&\tilde{b}_{21}=(0.13,-0.14,0.13,0.14,0.13,-0.14,0.13,0.14)^T;\\
&\tilde{b}_{22}=(-0.13,0.12,-0.13,0.12,-0.13,0.12,-0.13,0.12)^T;\\
&\tilde{I}_1=(-3,1,-3,1,-3,1,-3,1)^T; \tilde{I}_2=(2,4,2,4,2,4,2,4)^T.
\end{align*}

Moreover, for any $s=(s^0,s^1,\cdots,s^7)^T\in \mathbb{O}$, we define
\begin{align*}
&L_1(s)=\exp(-s^0-2s^1-s^2-2s^3-s^4-2s^5-s^6-2s^7),\\
&L_2(s)=\exp(-2s^0-s^1-2s^2-s^3-2s^4-s^5-2s^6-s^7),
\end{align*}
then
\begin{align*}
&f(s)\\
=&\frac{1-L_1(s)}{1+L_1(s)}e_0+\frac{1}{1+L_2(s)}e_1+\frac{1-L_2(s)}{1+L_2(s)}e_2+\frac{1}{1+L_1(s)}e_3\\
+&\frac{1-L_1(s)}{1+L_1(s)}e_4+\frac{1}{1+L_2(s)}e_5+\frac{1-L_2(s)}{1+L_2(s)}e_6+\frac{1}{1+L_1(s)}e_7\\
&g(s)\\
=&\frac{1}{1+L_2(s)}e_0+\frac{1-L_1(s)}{1+L_1(s)}e_1 +\frac{1}{1+L_1(s)}e_2+\frac{1-L_2(s)}{1+L_2(s)}e_3\\
+&\frac{1}{1+L_2(s)}e_4+\frac{1-L_1(s)}{1+L_1(s)}e_5 +\frac{1}{1+L_1(s)}e_6+\frac{1-L_2(s)}{1+L_2(s)}e_7		
\end{align*}

By some calculations, we get, for $q=1,2$
\begin{align*}
\overline{M}(q,\lambda)
=\mathrm{ones}(2,4)
%\left(
%\begin{array}{cccc}
%1&1&1&1\\
%1&1&1&1\end{array}\right)
\otimes
\left(
\begin{array}{llll}
0.5&0.5&1&0.25\\
1&0.25&0.5&0.5
\end{array}
\right)^T;\\
\overline{M}(q,\delta)
=\mathrm{ones}(2,4)
%\left(
%\begin{array}{cccc}
%1&1&1&1\\
%1&1&1&1\end{array}\right)
\otimes
\left(
\begin{array}{llll}
0.5&0.5&0.25&1\\
0.25&1&0.5&0.5
\end{array}
\right)^T.
\end{align*}
where $\mathrm{ones}(2,4)$ is a $2\times 4$ matrix with every element being $1$, and `$\otimes$' denotes the Kronecker product.

\begin{figure}
	\centering
	\includegraphics[width=0.45\textwidth]{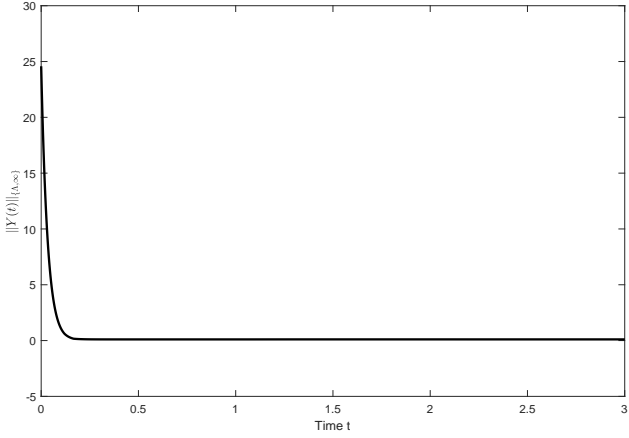}
	\caption{The trajectories of $\|Y(t)\|_{\{\Lambda,\infty\}}$ of OVNN (\ref{numerical-ovnn}) with constant time delays.}
	\label{fig1}
\end{figure}

Firstly, let us consider the case for constant time delay, $\tau_{11}=1, \tau_{12}=2, \tau_{21}=3, \tau_{22}=4$. Choose $\mu(t)=e^{0.02t}$, and $\Lambda=0.2\cdot\mathrm{ones}(16,1)$, so $\tau=4, \alpha=0.02$, calculations show
\begin{align*}
&\overline{T}_{\ell}(1,0.02,e^{0.08}-1)=(-0.8621,-1.2121,-0.9871,\\
&~~~~~~~-1.2221,-0.9671,-1.2421,-0.9521,-1.2421)<0,\\
&\overline{T}_{\ell}(2,0.02,e^{0.08}-1)=(-0.6781,-0.9281,-0.6881,\\
&~~~~~~~-0.9331,-0.6781,-0.8731,-0.6881,-0.9081)<0.
\end{align*}
As stated  by Corollary \ref{corollary-exponential stability}, OVNN (\ref{numerical-ovnn}) can be globally exponential-stable. Figure \ref{fig1} presents the trajectories of $\|Y(t)\|_{\{\Lambda,\infty\}}$ of system (\ref{numerical-ovnn}).
%, they all converge to 0 as time increases, which means there must exist a scalar $W_1>0$ such that
%	$$
%	\lVert Y(t)\rVert_{\{\xi,\infty\}}\leq \frac{W_1}{e^{0.02t}},
%	$$
%	i.e., the global exponential stability of system (\ref{numerical-ovnn}) is obtained.

Secondly, we consider unbounded time-varying delays. Define time delays as
\begin{align}\label{tau}
\tau_{11}(t)=\tau_{21}(t)=0.2t, \tau_{12}(t)=\tau_{22}(t)=0.1t
\end{align}
Choose $\mu(t)=t$, and $\Lambda=0.2\cdot\mathrm{ones}(16,1)$. In this case, $\alpha=0, \omega=0.2, \beta=0.25$, then calculations show
\begin{align*}
&\overline{T}_{\ell}(1,0,0.25)=(-0.59,-0.94,-0.715,\\
&~~~~~~~~~~~~~~~~-0.95,-0.695,-0.97,-0.68,-0.97)<0,\\
&\overline{T}_{\ell}(2,0,0.25)=(-0.37,-0.62,-0.38,\\
&~~~~~~~~~~~~~~~~-0.625,-0.37,-0.565,-0.38,-0.6)<0.
\end{align*}

\begin{figure}
	\centering
	\includegraphics[width=0.45\textwidth]{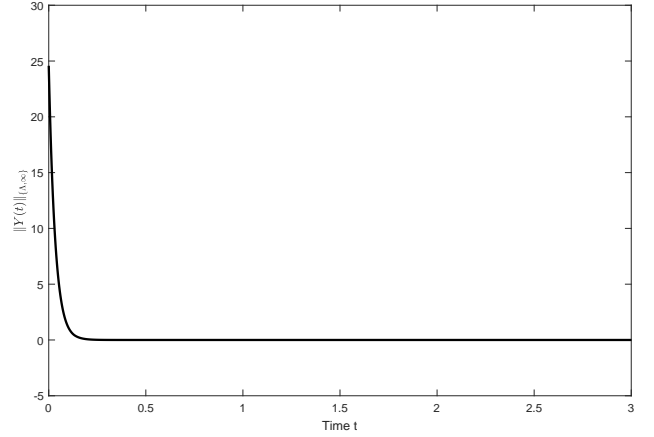}
	\caption{The trajectories of $\|Y(t)\|_{\{\Lambda,\infty\}}$ of OVNN (\ref{numerical-ovnn}) with unbounded time-varying delays.}
	\label{fig2}
\end{figure}

Based on Corollary \ref{corollary-power stability}, OVNN (\ref{numerical-ovnn}) has global power stability. Figure \ref{fig2} presents the trajectories of $\|Y(t)\|_{\{\Lambda,\infty\}}$ defined in (\ref{normyt}) of system (\ref{numerical-ovnn}).

\subsection{Finite-time control for OVNN with delay}
Consider a two-neuron OVNN system described as follows:
\begin{equation} \label{simu-ms}
\left\{
\begin{aligned}
\dot{w}_1(t)=&-d_1w_1(t)+a_{11}f_1(w_1(t)) \\
&+a_{12}f_2(w_2(t))+b_{11}g_1(w_1(t-\tau_{11}(t)))\\
&+b_{12}g_2(w_2(t-\tau_{12}(t)))+I_1+U_1(t), \\
\dot{w}_2(t)=&-d_2w_2(t)+a_{21}f_1(w_1(t)) \\
&+a_{22}f_2(w_2(t))+b_{21}g_1(w_1(t-\tau_{21}(t)))\\
&+b_{22}g_2(w_2(t-\tau_{22}(t)))+I_2+U_2(t),
\end{aligned}
\right.
\end{equation}
where $a_{ij}$ and $b_{ij}$ are defined as that in the above subsection, $\tau_{ij}$ are defined in (\ref{tau}), $i,j=1,2$. $d_1=0.1, d_2=0.2$ and
\begin{align*}
&\tilde{I}_1=(1.6,1.5,1.2,0.2,1.6,1.5,1.2,0.2)^T,\\
&\tilde{I}_2=(1.6,1.5,1.2,0.2,1.6,1.5,1.2,0.2)^T
\end{align*}
Moreover, for any $s=(s^0,s^1,\cdots,s^7)^T\in \mathbb{O}$, we define
\begin{align*}
&f(s)^{\ell}=\left\{\begin{array}{ll}
0.5\tanh(s^{\ell})+0.2\mathrm{sign}(s^{\ell}),&~~~\mathrm{even}~\ell\\
0.4\tanh(s^{\ell})-0.1\mathrm{sign}(s^{\ell}),&~~~\mathrm{odd}~\ell
\end{array}\right.\\
&g(s)^{\ell}=\left\{\begin{array}{ll}
0.35\tanh(s^{\ell})+0.05\mathrm{sign}(s^{\ell}),&\mathrm{even}~\ell\\
0.55\tanh(s^{\ell})-0.1\mathrm{sign}(s^{\ell}),&\mathrm{odd}~\ell
\end{array}\right.
\end{align*}

By some calculations, we get, for $q=1,2$
\begin{align*}
&\overline{M}(q,\lambda)=\mathrm{ones}(8,4)\otimes(0.7, 0.3);\\
&\overline{M}(q,\delta)=\mathrm{ones}(8,4)\otimes(0.4, 0.45).
\end{align*}

%	We define the initial values as
%	\begin{align*}
%	\tilde{w}_1(0)&=(-4.6,7,-7.6,2.4,-4.6,7,-7.6,2.4)^T,\\
%	\tilde{w}_2(0)&=(-2,  3,4.2, 5.1,-1.5,2,-4,  3.6)^T,
%	\end{align*}

\begin{figure}
	\centering
	\includegraphics[width=0.45\textwidth]{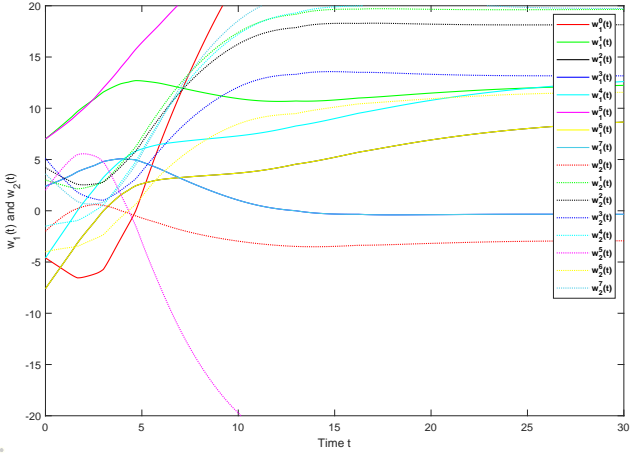}
	\caption{State trajectories of OVNN (\ref{simu-ms}) without control}
	\label{w1-w2 without control}
\end{figure}

Figure \ref{w1-w2 without control} presents the state trajectories of ${w}_p^{\ell}(t)$ for $p=1,2, \ell=0,1,\cdots,7$ of system (\ref{simu-ms}) without control, which means $U_p(t)=0$. We can easily find that the state trajectories of ${w}_p^{\ell}$ would not converge.

In the next, we add controller to make OVNN (\ref{simu-ms}) reach target points. We pick $\Lambda=\mathrm{ones}(4,1)\otimes(0.3,0.5, 0.5,0.3)^T$ and define $\mu(t)=t$. According to Assumption \ref{ass3}, $\alpha=0$, $\beta=0.25$.

At first, we choose the target as
$\tilde{\hat{Z}}^*=(0,\cdots,0)^T$, in accordance with (\ref{condition1})-(\ref{condition2}) in Theorem \ref{theorem3}, we obtain that
\begin{align*}
&\kappa_1^0>29.1733, \kappa_1^1>17.1880,  \kappa_1^2>28.4733, \kappa_1^3>17.1880,   \\
&\kappa_1^4>28.4733, \kappa_1^5>17.2180,  \kappa_1^6>28.4733, \kappa_1^7>17.1880,\\
&\kappa_2^0>17.3820, \kappa_2^1>29.9433  ,  \kappa_2^2>17.5220, \kappa_2^3> 29.9433 ,   \\
&\kappa_2^4>17.5220, \kappa_2^5> 30.0633 ,  \kappa_2^6>17.5220, \kappa_2^7>29.9433 ,\\
&\hat{\kappa}_1^0>4.0656   ,\hat{\kappa}_1^1> 3.9656,\hat{\kappa}_1^2>3.6656,\hat{\kappa}_1^3> 2.6656 ,\\
&\hat{\kappa}_1^4>4.0656  ,\hat{\kappa}_1^5>3.9656 ,\hat{\kappa}_1^6>3.6656  ,\hat{\kappa}_1^7>2.6656,\\
&\hat{\kappa}_2^0>4.4704  ,\hat{\kappa}_2^1> 4.3704  ,\hat{\kappa}_2^2>4.0704  ,\hat{\kappa}_2^3>3.0704,\\
&\hat{\kappa}_2^4>4.4704  ,\hat{\kappa}_2^5>4.3704 ,\hat{\kappa}_2^6>4.0704  ,\hat{\kappa}_2^7>3.0704.
\end{align*}
Therefore, the controllers $U_1(t),U_2(t)$ are designed as follows:
\begin{align}\label{simu-controller1}
U_1^{0}(t)&=-\mathrm{sign}(w_1^{0}(t))\big(29.2733|w_1^{0}(t)|+4.1656\big) \nonumber\\
U_1^{1}(t)&=-\mathrm{sign}(w_1^{1}(t))\big(17.2880|w_1^{1}(t)|+4.0656\big)\nonumber\\
U_1^{2}(t)&=-\mathrm{sign}(w_1^{2}(t))\big(28.5733|w_1^{2}(t)|+3.7656\big)\nonumber\\
U_1^{3}(t)&=-\mathrm{sign}(w_1^{3}(t))\big(17.2880|w_1^{3}(t)|+2.7656\big)\nonumber\\
U_1^{4}(t)&=-\mathrm{sign}(w_1^{4}(t))\big(28.5733|w_1^{4}(t)|+4.1656\big)\nonumber\\
U_1^{5}(t)&=-\mathrm{sign}(w_1^{5}(t))\big(17.3180|w_1^{5}(t)|+4.0656 \big)\nonumber\\
U_1^{6}(t)&=-\mathrm{sign}(w_1^{6}(t))\big(28.5733|w_1^{6}(t)|+3.7656 \big)\nonumber\\
U_1^{7}(t)&=-\mathrm{sign}(w_1^{7}(t))\big(17.2880|w_1^{7}(t)|+2.7656\big) \nonumber\\
U_2^{0}(t)&=-\mathrm{sign}(w_2^{0}(t))\big(17.4820|w_2^{0}(t)|+4.5704\big) \nonumber\\
U_2^{1}(t)&=-\mathrm{sign}(w_2^{1}(t))\big(30.0433|w_2^{1}(t)|+4.4704\big)\nonumber\\
U_2^{2}(t)&=-\mathrm{sign}(w_2^{2}(t))\big(17.6220|w_2^{2}(t)|+4.1704\big)\nonumber\\
U_2^{3}(t)&=-\mathrm{sign}(w_2^{3}(t))\big(30.0433|w_2^{3}(t)|+3.1704\big)\nonumber\\
U_2^{4}(t)&=-\mathrm{sign}(w_2^{4}(t))\big(18.6220|w_2^{4}(t)|+4.5704\big)\nonumber\\
U_2^{5}(t)&=-\mathrm{sign}(w_2^{5}(t))\big(30.1633|w_2^{5}(t)|+4.4704\big)\nonumber\\
U_2^{6}(t)&=-\mathrm{sign}(w_2^{6}(t))\big(17.6220|w_2^{6}(t)|+4.1704\big)\nonumber\\
U_2^{7}(t)&=-\mathrm{sign}(w_2^{7}(t))\big(30.0433|w_2^{7}(t)|+3.1704\big)\nonumber\\
\end{align}

\begin{figure}
	\centering
	\includegraphics[width=0.45\textwidth]{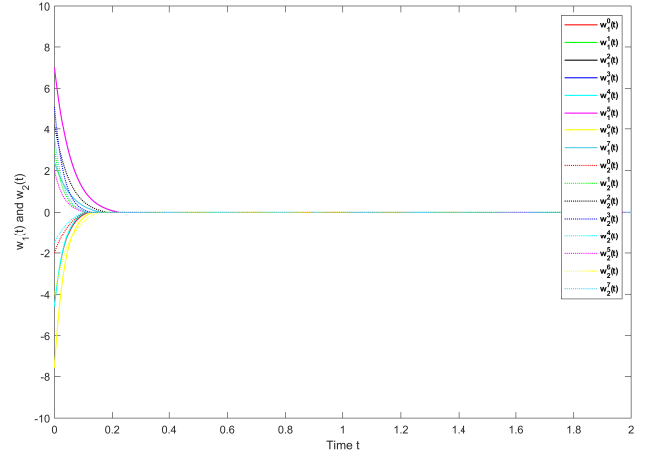}
	\caption{State trajectories of OVNN (\ref{simu-ms}) under control (\ref{simu-controller1})}
	\label{w1-w2 under control1}
\end{figure}

\begin{figure}
	\centering
	\includegraphics[width=0.45\textwidth]{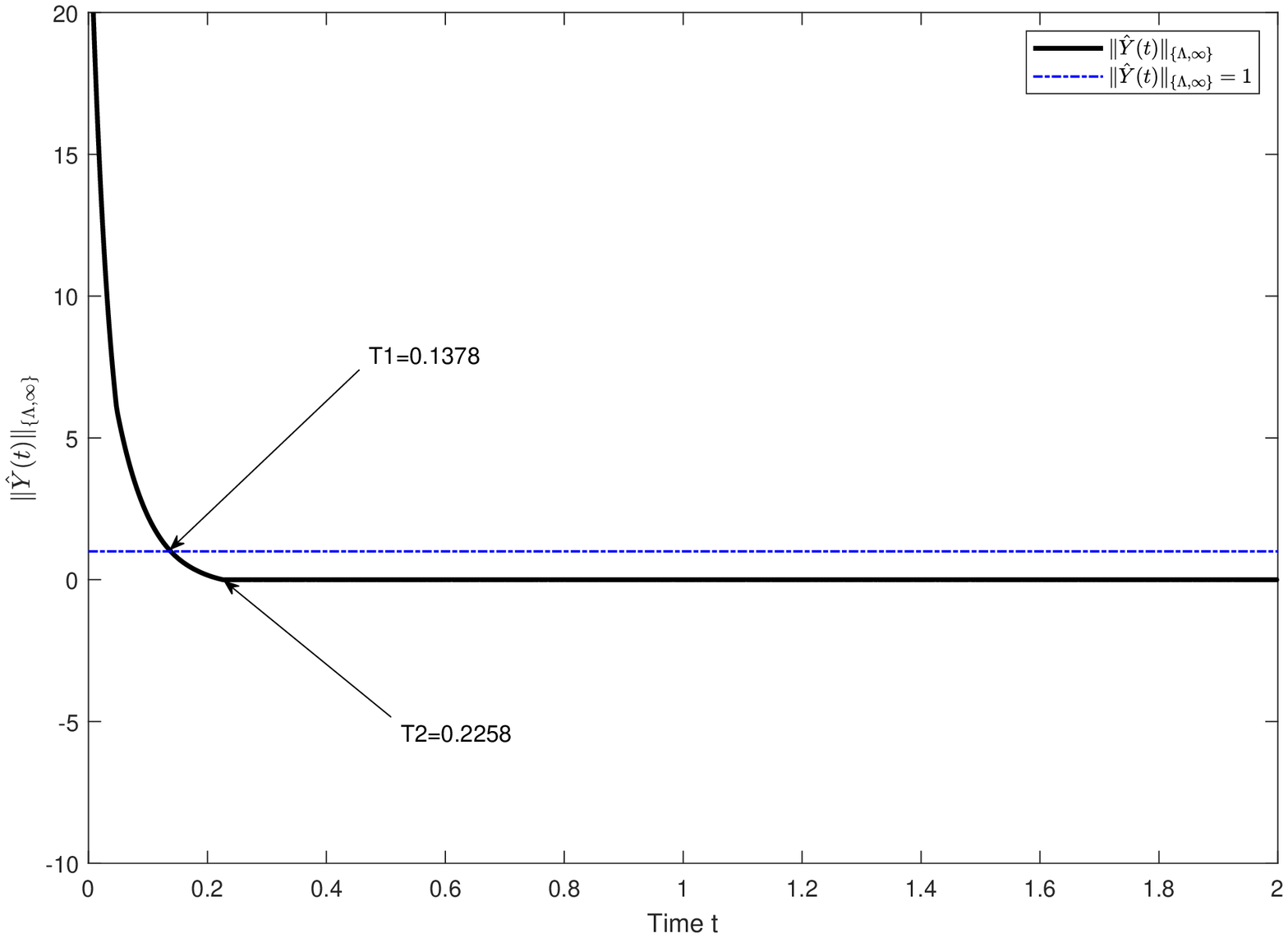}
	\caption{Trajectory of $\lVert \hat{Y}(t)\rVert_{\{\Lambda,\infty\}}$ of OVNN (\ref{simu-ms}) under control (\ref{simu-controller1})}
	\label{yhat-norm under control1}
\end{figure}

Figures \ref{w1-w2 under control1} and \ref{yhat-norm under control1} present the state trajectories and $\lVert \hat{Y}(t)\rVert_{\{\Lambda,\infty\}}$ defined in (\ref{normhyt}) of system (\ref{simu-ms}) under control (\ref{simu-controller1}). It is obvious that the state trajectories of $w_p^{\ell}$ for $p=1,2, \ell=0,1,\cdots,7$ converge to the target as $t \rightarrow \infty$. Consistent with what we have proved, $\lVert \hat{Y}(t)\rVert_{\{\Lambda,\infty\}}$ firstly decreases from initial state to 1 in finite time $T_1$, then decreases to 0 before $T_2$, i.e., OVNNs (\ref{simu-ms}) reach stability in finite time.

\begin{figure}
	\centering
	\includegraphics[width=0.45\textwidth]{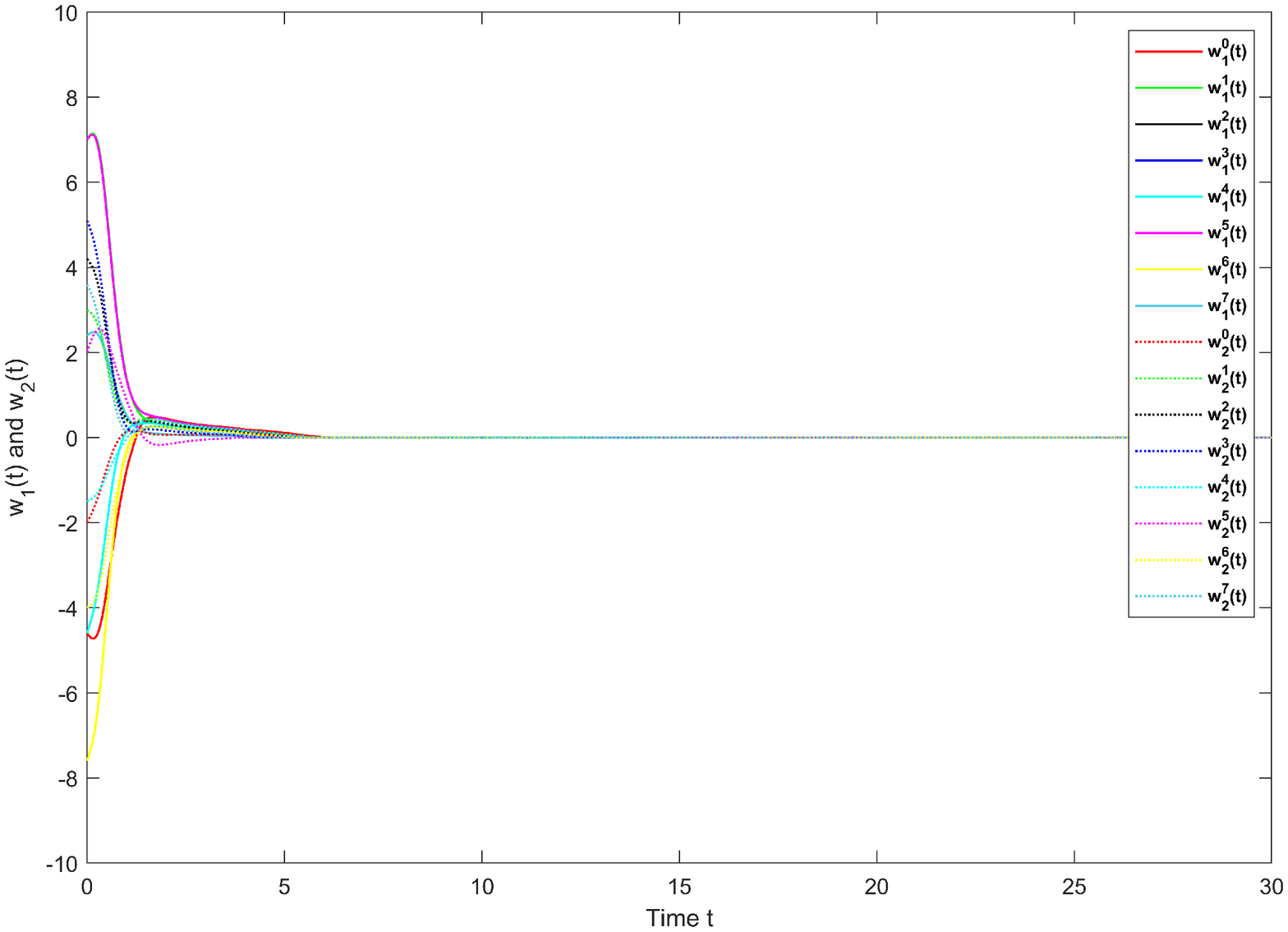}
	\caption{State trajectories of OVNN (\ref{simu-ms}) under adaptive control (\ref{add10})-(\ref{simu-rules of kappa})}
	\label{w1-w2 under control2}
\end{figure}

\begin{figure}
	\centering
	\includegraphics[width=0.45\textwidth]{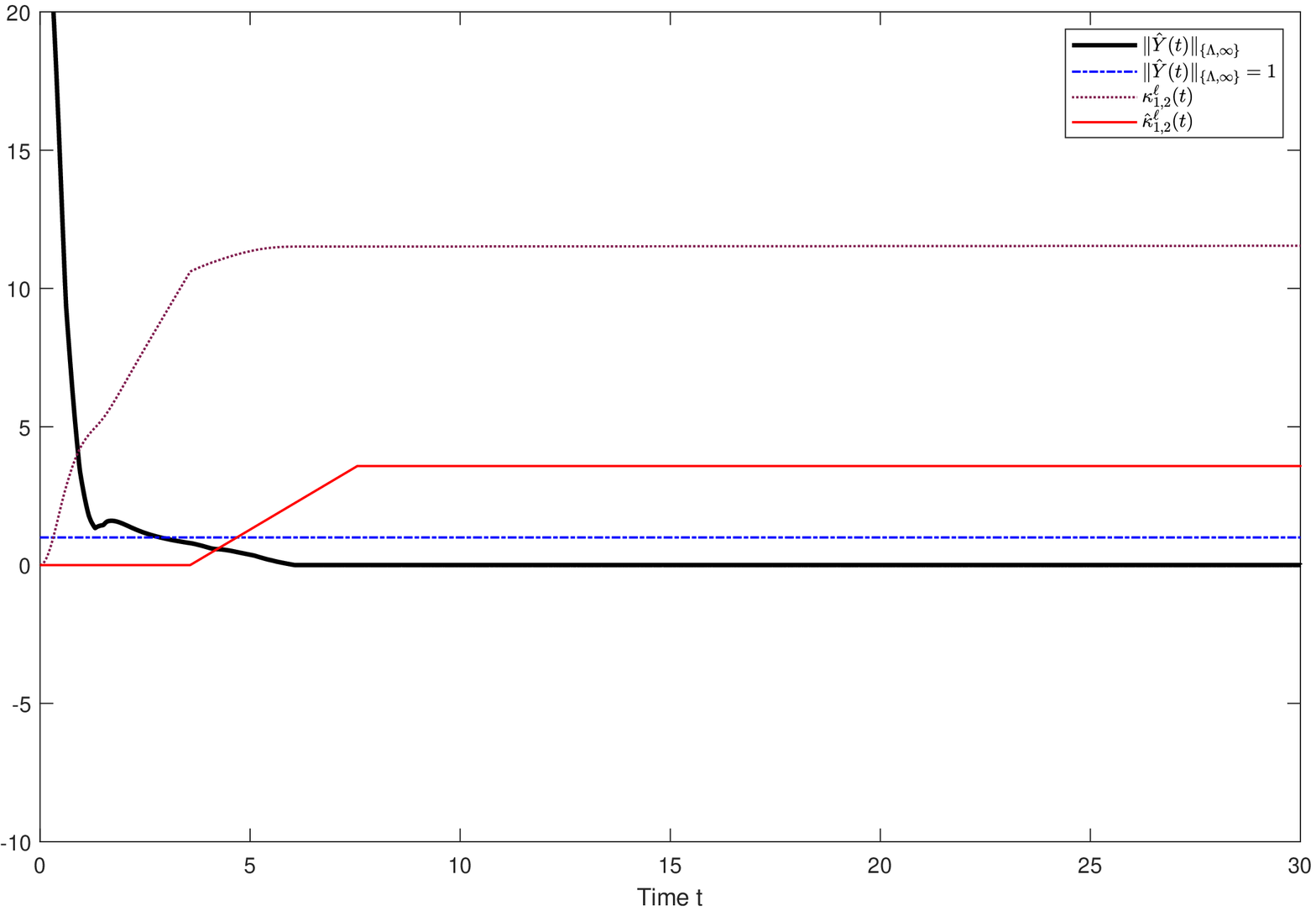}
	\caption{Trajectory of $\lVert \hat{Y}(t)\rVert_{\{\Lambda,\infty\}}$, $\hat{\kappa}(t)$ and $\kappa(t)$ of OVNN (\ref{simu-ms}) under adaptive control (\ref{add10})-(\ref{simu-rules of kappa})}
	\label{yhat-norm under control2}
\end{figure}

Then we apply the adaptive strategy on (\ref{simu-ms}). The controllers $U_1$ and $U_2$ are defined as
\begin{align}\label{add10}
U_p^{\ell}(t)=-\mathrm{sign}(\hat{w}_p^{\ell}(t))\big(\kappa(t)|\hat{w}_p^{\ell}(t)|+\hat{\kappa}(t)\big),~~p=1,2
\end{align}
with adaptive rules
\begin{align}\label{simu-rules of chi}
\dot{\hat{\kappa}}(t)=
\left \{
\begin{array}{ll}
0, &\sup_{t-0.2t\leq s \leq t}\lVert \hat{Y}(s)\rVert_{\{\Lambda,\infty\}} >1\\
0.9, &\sup_{t-0.2t\leq s \leq t}\lVert \hat{Y}(s)\rVert_{\{\Lambda,\infty\}} \leq 1 \\
0,  &\sup_{t-0.2t\leq s \leq t}\lVert \hat{Y}(s)\rVert_{\{\Lambda,\infty\}} = 0
\end{array}
\right .
\end{align}
and
\begin{align}\label{simu-rules of kappa}
\dot{\kappa}(t)=
\left \{
\begin{array}{ll}
0.9t\lVert \hat{Y}(t)\rVert_{\{\Lambda,\infty\}},\sup_{t-0.2t\leq s \leq t}\lVert \hat{Y}(s)\rVert_{\{\Lambda,\infty\}} >1\\
0.9\lVert \hat{Y}(t)\rVert_{\{\Lambda,\infty\}},~\sup_{t-0.2t\leq s \leq t}\lVert \hat{Y}(s)\rVert_{\{\Lambda,\infty\}} \leq 1
\end{array}
\right .
\end{align}
where $\hat{\kappa}(0)=\kappa(0)=0$. According to Theorem \ref{theorem4}, the adaptive finite time stablity of (\ref{simu-ms}) can be realized. Figure \ref{w1-w2 under control2}  presents the state trajectories of system (\ref{simu-ms}) under adaptive controllers (\ref{add10})-(\ref{simu-rules of kappa}). Figure \ref{yhat-norm under control2} illustrates the $\lVert \hat{Y}(t)\rVert_{\{\Lambda,\infty\}}$, $\hat{\kappa}(t)$ and $\kappa(t)$ of system (\ref{simu-ms}) under adaptive controllers (\ref{add10})-(\ref{simu-rules of kappa}).

To show that our algorithm can drive OVNN to any given point, in the next, we choose another point, $\tilde{\hat{Z}}^*=(1,2,\cdots,16)^T$. With the same calculation in the above, controllers $U_1(t), U_2(t)$ are designed as
\begin{align}\label{simu-controller1-ep2}
U_1^{0}(t)&=-\mathrm{sign}(w_1^{0}(t))\big(29.2733|w_1^{0}(t)|+6.4788\big) \nonumber\\
U_1^{1}(t)&=-\mathrm{sign}(w_1^{1}(t))\big(17.2880|w_1^{1}(t)|+5.2113\big)\nonumber\\
U_1^{2}(t)&=-\mathrm{sign}(w_1^{2}(t))\big(28.5733|w_1^{2}(t)|+5.0113\big)\nonumber\\
U_1^{3}(t)&=-\mathrm{sign}(w_1^{3}(t))\big(17.2880|w_1^{3}(t)|+4.1113\big)\nonumber\\
U_1^{4}(t)&=-\mathrm{sign}(w_1^{4}(t))\big(28.5733|w_1^{4}(t)|+5.6113\big)\nonumber\\
U_1^{5}(t)&=-\mathrm{sign}(w_1^{5}(t))\big(17.3180|w_1^{5}(t)|+8.6244 \big)\nonumber\\
U_1^{6}(t)&=-\mathrm{sign}(w_1^{6}(t))\big(28.5733|w_1^{6}(t)|+5.4113 \big)\nonumber\\
U_1^{7}(t)&=-\mathrm{sign}(w_1^{7}(t))\big(17.2880|w_1^{7}(t)|+4.5113\big) \nonumber\\
U_2^{0}(t)&=-\mathrm{sign}(w_2^{0}(t))\big(17.4820|w_2^{0}(t)|+3.5750\big) \nonumber\\
U_2^{1}(t)&=-\mathrm{sign}(w_2^{1}(t))\big(30.0433|w_2^{1}(t)|+8.7339\big)\nonumber\\
U_2^{2}(t)&=-\mathrm{sign}(w_2^{2}(t))\big(17.6220|w_2^{2}(t)|+8.6339\big)\nonumber\\
U_2^{3}(t)&=-\mathrm{sign}(w_2^{3}(t))\big(30.0433|w_2^{3}(t)|+7.8339\big)\nonumber\\
U_2^{4}(t)&=-\mathrm{sign}(w_2^{4}(t))\big(18.6220|w_2^{4}(t)|+9.4339\big)\nonumber\\
U_2^{5}(t)&=-\mathrm{sign}(w_2^{5}(t))\big(30.1633|w_2^{5}(t)|+6.7698\big)\nonumber\\
U_2^{6}(t)&=-\mathrm{sign}(w_2^{6}(t))\big(17.6220|w_2^{6}(t)|+7.5889\big)\nonumber\\
U_2^{7}(t)&=-\mathrm{sign}(w_2^{7}(t))\big(30.0433|w_2^{7}(t)|+10.4789\big)\nonumber\\
\end{align}
\begin{figure}
	\centering
	\includegraphics[width=0.45\textwidth]{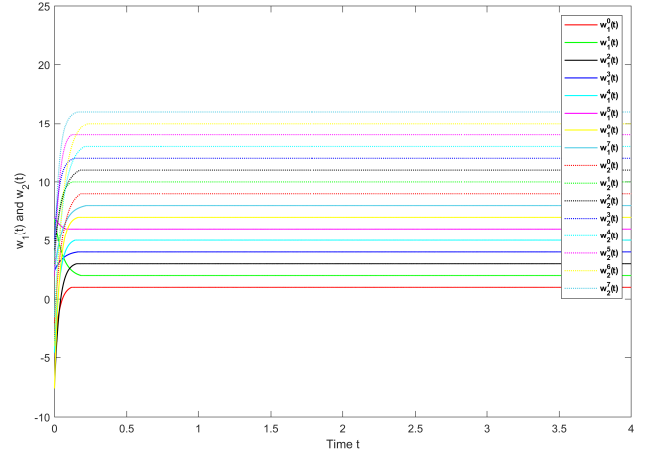}
	\caption{State trajectories of OVNN (\ref{simu-ms}) under control (\ref{simu-controller1-ep2})}
	\label{w1-w2 under control1-ep2}
\end{figure}
\begin{figure}
	\centering
	\includegraphics[width=0.45\textwidth]{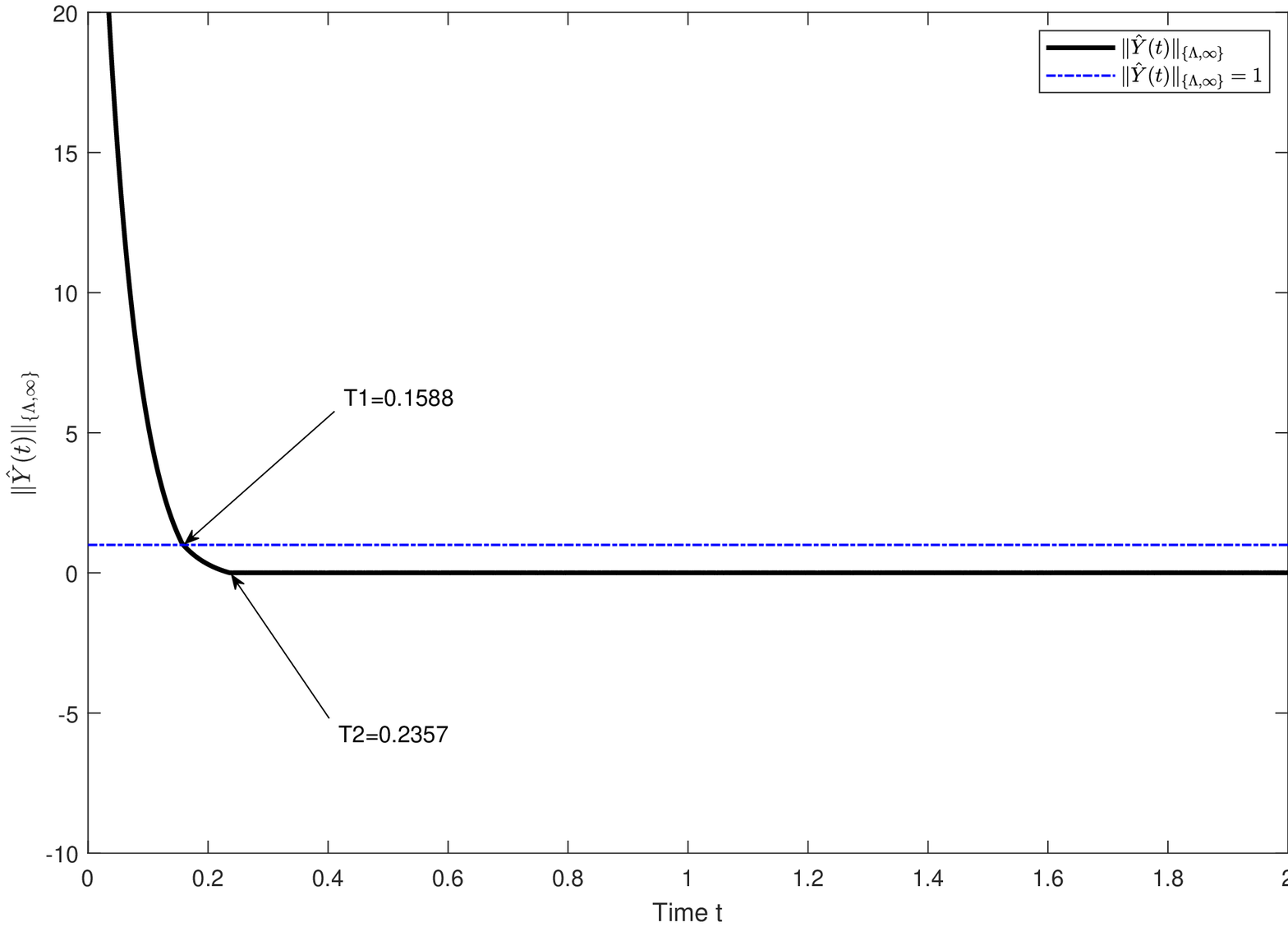}
	\caption{Trajectory of $\lVert \hat{Y}(t)\rVert_{\{\Lambda,\infty\}}$ of OVNN (\ref{simu-ms}) under control (\ref{simu-controller1-ep2})}
	\label{yhat-norm under control1-ep2}
\end{figure}

Figure \ref{w1-w2 under control1-ep2} shows the state trajectories of system (\ref{simu-ms}) under control (\ref{simu-controller1-ep2}). The state trajectories of $w_p^{\ell}$ for $p=1,2, \ell=0,1,\cdots,7$ converge to the target points. Figure \ref{yhat-norm under control1-ep2} shows  the trajectories of $\lVert \hat{Y}(t)\rVert_{\{\Lambda,\infty\}}$, which decreases from initial state to $1$, then to $0$ in finite time.

\begin{figure}
	\centering
	\includegraphics[width=0.45\textwidth]{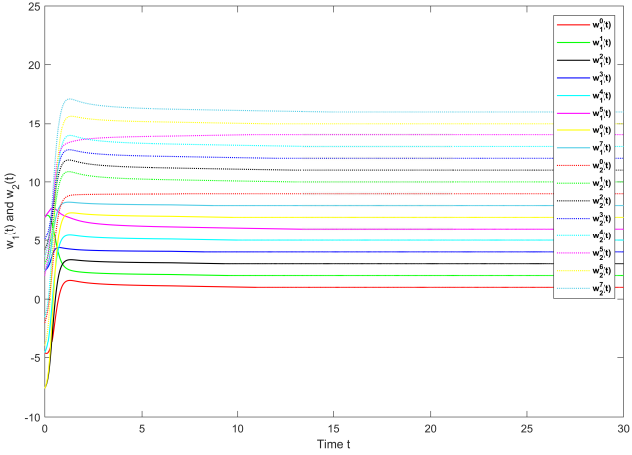}
	\caption{State trajectories of OVNN (\ref{simu-ms}) under adaptive control (\ref{add10})-(\ref{simu-rules of kappa})}
	\label{w1-w2 under control2-ep2}
\end{figure}

\begin{figure}
	\centering
	\includegraphics[width=0.45\textwidth]{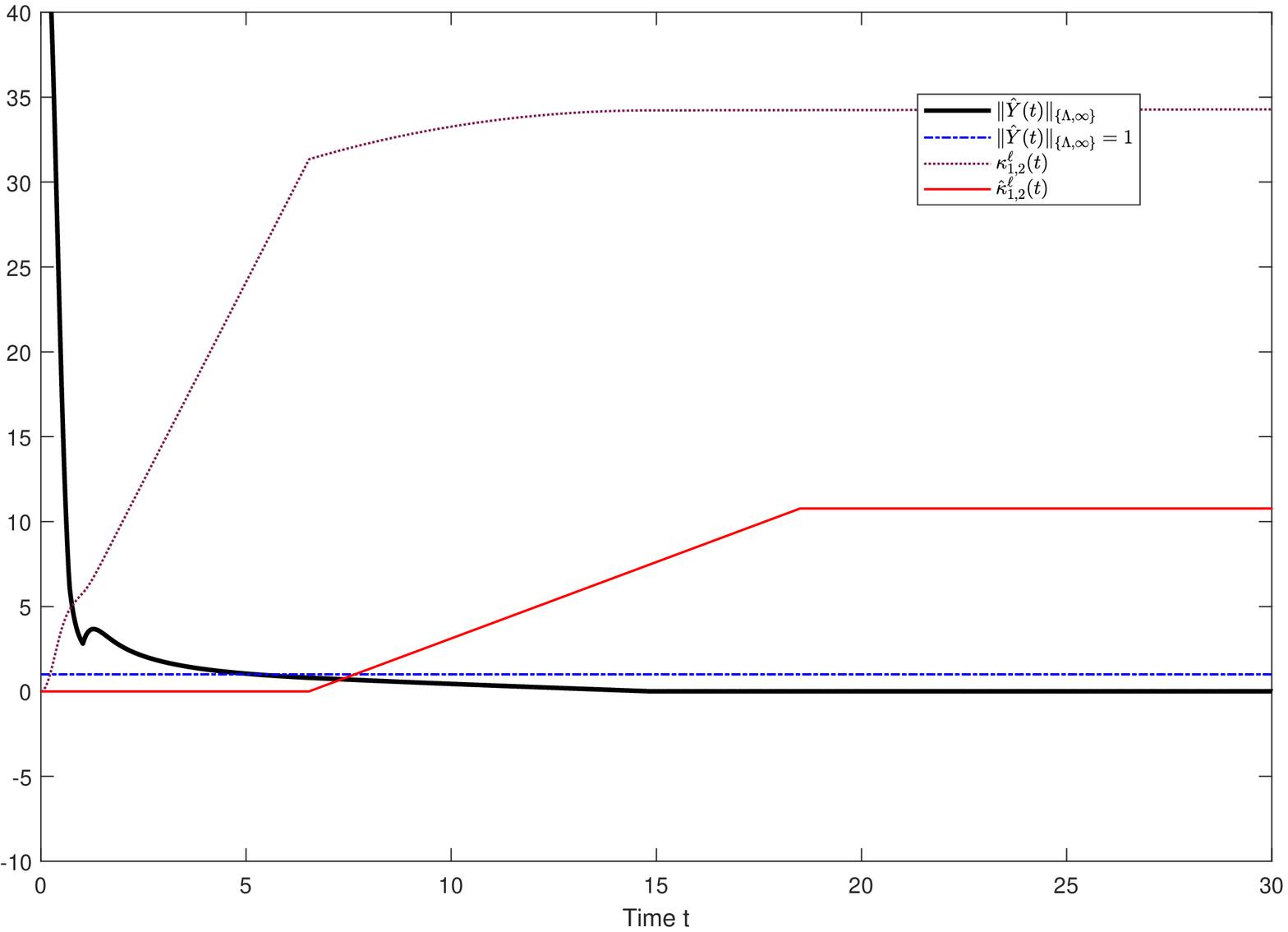}
	\caption{Trajectory of $\lVert \hat{Y}(t)\rVert_{\{\Lambda,\infty\}}$, $\hat{\kappa}(t)$ and $\kappa(t)$ of OVNN (\ref{simu-ms}) under adaptive control (\ref{add10})-(\ref{simu-rules of kappa})}
	\label{yhat-norm under control2-ep2}
\end{figure}

Moreover, for adaptive control (\ref{add10})-(\ref{simu-rules of kappa})
can also be used. Figure \ref{w1-w2 under control2-ep2} shows the state trajectories of OVNN (\ref{simu-ms}) with adaptive  rules (\ref{add10})-(\ref{simu-rules of kappa}). It is clear that OVNN (\ref{simu-ms}) reaches the second target in finite time.  With this target,  Figure \ref{yhat-norm under control2-ep2} illustrates the $\lVert \hat{Y}(t)\rVert_{\{\Lambda,\infty\}}$, $\hat{\kappa}(t)$ and $\kappa(t)$ of system (\ref{simu-ms}) under adaptive control (\ref{add10})-(\ref{simu-rules of kappa}).	

\section{Conclusion}\label{con}
	In this paper, we discuss the issues of  the $\mu$-stability and finite-time stability of OVNNs with unbounded and asynchronous time-varying delays by using decomposing skill and generalized norm. Some criteria are established for the $\mu$-stability of the delayed OVNN. Moreover, once the equilibrium is not stable, we propose controllers only depending on the system state, we also design the adaptive technique to realized finite-time stability. The validity of proposed algorithm is proved by 2PM. Finally, the effectiveness of the obtained criteria is shown by numerical examples. Recently, a generalization of both the complex and quaternion algebras is the Clifford algebra and Clifford-valued neural networks (ClVNNs) have attracted researchers' interests. However, Clifford algebra does not have the property of being a normed division algebra, i.e., a norm and a multiplicative inverse cannot be defined for them, see \cite{popa2}. Thus the obtained results need further research for ClVNNs with unbounded time-varying delays.

\section*{Appendix A: Proof of Theorem \ref{theorem1}}
Let us discuss the companion network of OVNN (\ref{equation5})
\begin{align}\label{conetwork}
\dot{z}_p(t)=&-d_pz_p(t)+\sum_{q=1}^na_{pq}f_q(z_q(t))+\sum_{q=1}^nb_{pq}g_q(z_q(t)))\nonumber\\
&+I_p,~~p=1,\cdots,n
\end{align}

It is clear that the equilibrium point of OVNN (\ref{equation5}) is the same as that of (\ref{conetwork}). If inequalities (\ref{t0}) hold, there is a constant $\epsilon>0$ which is small enough to make the following inequalities hold:
\begin{align*}
		S_{\ell}(p)=\epsilon\Lambda_p[\ell]+T_{\ell}(p)\leq0, ~~{\ell}=0,1,\cdots,7
\end{align*}
		
Similarly, NN (\ref{conetwork}) can be decomposed into eight RVNNs
\begin{align}\label{decomposed-conetwork}
\dot{z}_p^{\ell}(t)=-d_pz_p^{\ell}+\sum\limits_{q=1}^n(\tilde{a}_{pq}^TM^{\ell}\tilde{f}_q+\tilde{b}_{pq}^TM^{\ell}\tilde{g}_q)+I_p^{\ell},
\end{align}
		
Define $u_p^{\ell}(t)= \dot{z}_p^{\ell}(t)$, $\ell=0,1,\cdots,7$, and obviously $\dot{\widetilde{z}}_p(t)={\widetilde{u}}_p(t)$. In accordance with Lemma \ref{lemma2}, we get
\begin{align*}
		\dot{u}_p^{\ell}(t)=&-d_pu_p^{\ell}(t)\\	&+\sum\limits_{q=1}^n\Big[\tilde{a}_{pq}^TM^{\ell}\overline{M}(f_q)+\tilde{b}_{pq}^TM^{\ell}\overline{M}(g_q)\Big]{\widetilde{u}}_p(t)
		\end{align*}
		
Define
		\begin{align*}
		Z(t)=&(z_1^0,\cdots,z_n^0,z_1^1,\cdots,z_n^1,z_1^2,\cdots,z_n^2,z_1^3,\cdots,z_n^3,\\
		&~z_1^4,\cdots,z_n^4,z_1^5,\cdots,z_n^5,z_1^6,\cdots,z_n^6,z_1^7,\cdots,z_n^7)^T\\
		X(t)=&(u_1^0,\cdots,u_n^0,u_1^1,\cdots,u_n^1,u_1^2,\cdots,u_n^2,u_1^3,\cdots,u_n^3,\\
		&~u_1^4,\cdots,u_n^4,u_1^5,\cdots,u_n^5,u_1^6,\cdots,u_n^6,u_1^7,\cdots,u_n^7)^T
		\end{align*}
and
		\begin{align*}
		\lVert X(t)\rVert_{\{\Lambda,\infty\}}=\max_{\ell}\Big\{
		&\max_p\{\Lambda_p[\ell]^{-1}|u_p^{\ell}(t)|\}\Big\}
		\end{align*}
		
For time $t$, without loss of generality, suppose $p_0=p_0(t)$ is an index which satisfies $\lVert X(t)\rVert_{\{\Lambda,\infty\}}=|\xi_{p_0}^{-1}u_{p_0}^0|$, then
		\begin{align*}
		&\xi_{p_0}\frac{d\lVert e^{\epsilon t}X(t)\rVert_{\{\Lambda,\infty\}}}{dt}=\frac{d|e^{\epsilon t}u_{p_0}^0(t)|}{dt}\\
%		=&\epsilon e^{\epsilon t}|u_{p_0}^0(t)|+e^{\epsilon t}\mathrm{sign}\{u_{p_0}^0(t)\} \times \Big\{-d_{p_0}u_{p_0}^0(t)\\ &+\sum\limits_{q=1}^n\Big[\tilde{a}_{p_0q}^TM^0\overline{M}(f_q)+\tilde{b}_{p_0q}^TM^0\overline{M}(g_q)\Big]{\widetilde{u}}_q(t)		\Big\}\\
		=&\epsilon e^{\epsilon t}|u_{p_0}^0(t)|+e^{\epsilon t}\mathrm{sign}\{u_{p_0}^0(t)\} \times \Big\{-d_{p_0}u_{p_0}^0(t)\\ &+\Big[\tilde{a}_{p_0p_0}^TM^0\overline{M}(f_{p_0})+\tilde{b}_{p_0p_0}^TM^0\overline{M}(g_{p_0})\Big]{\widetilde{u}}_{p_0}(t)\\
		&+\sum\limits_{q\neq p_0}\Big[\tilde{a}_{p_0q}^TM^0\overline{M}(f_q)+\tilde{b}_{p_0q}^TM^0\overline{M}(g_q)\Big]{\widetilde{u}}_q(t)		
		\Big\}\\
		=&\epsilon e^{\epsilon t}|u_{p_0}^0(t)|+e^{\epsilon t}\mathrm{sign}\{u_{p_0}^0(t)\} \times \Big\{-d_{p_0}u_{p_0}^0(t)\\		&+\Big[\tilde{a}_{p_0p_0}^TM^0\overline{M}(f_{p_0},0)+\tilde{b}_{p_0p_0}^TM^0\overline{M}(g_{p_0},0)\Big]{\widetilde{u}}_{p_0}(t)\\
		&+\tilde{a}_{p_0p_0}^TM^0\sum\limits_{{\ell}=1}^7\overline{M}(f_{p_0},\ell){\widetilde{u}}_{p_0}(t)\\
		&+\tilde{b}_{p_0p_0}^TM^0\sum\limits_{\ell=1}^7\overline{M}(g_{p_0},\ell){\widetilde{u}}_{p_0}(t)\\
		&+\sum\limits_{q\neq p_0}\Big[\tilde{a}_{p_0q}^TM^0\overline{M}(f_q)+\tilde{b}_{p_0q}^TM^0\overline{M}(g_q)\Big]{\widetilde{u}}_q(t)		
		\Big\}\\
		\leq&\epsilon e^{\epsilon t}|u_{p_0}^0(t)|+e^{\epsilon t}\Big\{\Big[-d_{p_0}+\{\tilde a_{p_0p_0}^T  M^0\}^+\overline{M}_0(p_0,\lambda) \\
		&+\{\tilde b_{p_0p_0}^T M^0\}^+\overline{M}_0(p_0,\delta)\Big]|u_{p_0}^0| \\
		%2+3
		&+\Big[|\tilde a_{p_0p_0}^T M^0|\overline{M}(p_0,\lambda)+|\tilde b_{p_0p_0}^T M^0|\overline{M}(p_0,\delta)\Big]|u_{p_0}^{[0]}|\\
		&+\sum\limits_{q\neq p_0}\Big[|\tilde a_{p_0q}^T M^0|\overline{M}(q,\lambda)+|\tilde b_{p_0q}^T M^0|\overline{M}(q,\delta)\Big]|u_q|\Big\}\\
\le&S_0(p_0)\lVert e^{\epsilon t}X(t) \rVert_{\{\Lambda,\infty\}} \leq 0.
		\end{align*}

		Accordingly, $\lVert X(t)\rVert_{\{\Lambda,\infty\}}=\lVert \dot{Z}(t)\rVert_{\{\Lambda,\infty\}}=\mathit{O}(e^{-\epsilon t})$. So for time $t_1>t_2$, there exists a constant $W>0$ such that
		\begin{align*}
		\lVert Z(t_1)&-Z(t_2)\rVert_{\{\Lambda,\infty\}}=\bigg\lVert \int_{t_2}^{t_1} X(t) dt\bigg\rVert_{\{\Lambda,\infty\}}\\
		&\leq \int_{t_2}^{t_1} We^{-\epsilon t}\,dt =\frac{W}{\epsilon}(e^{-\epsilon t_2}-e^{-\epsilon t_1}) \leq \frac{W}{\epsilon}e^{-\epsilon t_2}
		\end{align*}
		
According to Cauchy convergence principle, we obtain $\lim_{t \rightarrow \infty}Z(t)= Z^*$, where $Z^*$ is an equilibrium of OVNN (\ref{conetwork}), also the equilibrium for OVNN (\ref{equation5}).
		
		Now, by replacing the parameters of the above inequality, we can prove that
		$$
		\lVert Z(t)-Z^*\rVert_{\{\Lambda,\infty\}}=\bigg\lVert \int_{t}^{\infty} \dot{Z}(t) dt \bigg\rVert_{\{\Lambda,\infty\}}\leq \frac{W}{\epsilon}e^{-\epsilon t}
		$$
so we prove the uniqueness of the equilibrium point	above and obtain that any $Z(t)$ exponentially converges to ${Z^*}$.
		
\section*{Appendix B: Proof of Theorem \ref{theorem2}}
When inequalities (\ref{tt0}) hold, we can find that inequalities (\ref{t0}) hold. Then according to Theorem \ref{theorem1},
OVNN (\ref{equation5}) can reach a unique equilibrium point $Z^*=(z_1^*,z_2^*,\cdots,z_n^*)\in \mathbb{O}^n$.
		
		According to (\ref{ab-equation}) and (\ref{tt0}), there exists a sufficiently large $t^\star>0 $, such that for any $t>t^\star$,
		\begin{align}
		\overline{T}_\ell\bigg(p,\frac{\dot{\mu}(t)}{\mu(t)}, \frac{\mu(t)}{\mu(t-\tau(t))}-1 \bigg)<0,~~ \ell=0,1,\cdots,7
		\end{align}
		
		Define
		\begin{align}\label{useful3}
		P(t)=\sup_{t-\tau(t)<s\leq t}\lVert\mu(s)Y(s) \rVert_{\{\Lambda,\infty\}}
		\end{align}
		
For $t>t^{\star}$, if $\lVert \mu(t)Y(t)\rVert_{\{\Lambda,\infty\}}<P(t)$, then we can find an interval $\delta$ to make  $\lVert\mu(s)Y(s)\rVert_{\{\Lambda,\infty\}}<P(t)$ hold when $s \in (t,t+\delta)$. Otherwise, if there exists a time point $t_0$ such that $\lVert \mu(t_0)Y(t_0)\rVert_{\{\Lambda,\infty\}}=P(t_0)$, without loss of generality, we suppose $p_0=p_0(t_0)$ is an index which satisfies $\lVert Y(t_0)\rVert_{\{\Lambda,\infty\}}=|\xi_{p_0}^{-1}\overline{w}_{p_0}^0(t_0)|$, we obtain
		\begin{align*}
		&\xi_{p_0}\frac{d\lVert \mu(t)Y(t)\rVert_{\{\Lambda,\infty\}}}{dt}\bigg|_{t=t_0}
		=\frac{d|\mu(t)\overline{w}_{p_0}^0(t)|}{dt}\bigg|_{t=t_0}\\
		=&\dot{\mu}(t_0)|\overline{w}_{p_0}^0(t_0)|+\mu(t_0)\mathrm{sign}\{\overline{w}_{p_0}^0(t_0)\} \times \bigg\{-d_{p_0}\overline{w}_{p_0}^0(t_0)\\
		&+\tilde{a}_{p_0p_0}^TM^0(\tilde f_{p_0}-\tilde f_{p_0}^*)+\sum\limits_{q \neq p_0}\tilde{a}_{p_0q}^TM^0(\tilde f_q-\tilde f_q^*)\\
		&+\sum\limits_{q=1}^n\tilde{b}_{p_0q}^TM^0(\tilde g_{q\tau_{pq}}-\tilde g_q^*)\bigg\}\\
		%小于等于
		%a_p1p1
		\leq&\mu(t_0)\bigg\{\xi_{p_0} \bigg[-d_{p_0}+\frac{\dot{\mu}(t_0)}{\mu(t_0)}+\{\tilde a_{p_0p_0}^T M^0\}^+\overline{M}_0(p_0,\lambda)\bigg]\\
		&\cdot|\xi_{p_0}^{-1}\overline{w}_{p_0}^0(t_0)|+|\tilde a_{p_0p_0}^T M^0|\overline{M}(p_0,\lambda)\mathrm{diag}(\Lambda_{p_0}^{[0]})|\overline{Y}_p(t_0)^{[0]}|\\
		&+\sum\limits_{q\neq p_0}|\tilde a_{p_0q}^T M^0|\overline{M}(q,\lambda)\mathrm{diag}(\Lambda_q)|\overline{Y}_q(t_0)| \bigg\}\\
		&+\mu(t_0-\tau_{p_0q}(t_0))\frac{\mu(t_0)}{\mu(t_0-\tau_{p_0q}(t_0))}\sum\limits_{q=1}^n(|\tilde b_{p_0q}^T M^0|)\\
		&~~~\cdot\overline{M}(q,\delta)\mathrm{diag}(\Lambda_q)|\overline{Y}_q(t_0-\tau_{p_0q}(t_0))|\\
\leq&\overline{T}_0\bigg(p_0,\frac{\dot{\mu}(t_0)}{\mu(t_0)}, \frac{\mu(t_0)}{\mu(t_0-\tau(t_0))}-1\bigg)P(t_0)<0.
		\end{align*}
where $\overline{Y}_p(t)=(\Lambda_p[0]^{-1}\overline{w}_p^{0}(t),\cdots,\Lambda_p[7]^{-1}\overline{w}_p^{7}(t))^T$.

In summary, we conclude that $P(t)=P(t^\star)$ for all $t>t^\star$, which implies $\lVert \mu(t)Y(t)\rVert_{\{\Lambda,\infty\}}\leq P(t^\star)$ and
		$$
		\lVert Y(t)\rVert_{\{\Lambda,\infty\}}\leq \frac{P(t^\star)}{\mu(t)},
		$$

Therefore, OVNN (\ref{equation5}) is globally $\mu$-stable.
		
\section*{Appendix C: Proof of Theorem \ref{theorem3}}
	According to Assumption \ref{ass3} and (\ref{condition1}), we can find a time $T$, such that for $t\ge T$,
	\begin{align}\label{proof-condition1}
\kappa_p^{\ell}>\Lambda_p[\ell]^{-1}\overline{T}_{\ell}\bigg(p,\frac{\dot{\mu}(t)}{\mu(t)}, \frac{\mu(t)}{\mu(t-\tau(t))}-1\bigg)
	\end{align}
In the following, we always discuss from $T$. The proof process is divided into two phases.
	
\textbf{Phase I}: The norm $\lVert \hat{Y}(t)\rVert_{\{\Lambda,\infty\}}$ defined in (\ref{normhyt}) will be proved to evolve from initial values to 1 in finite time $T_1$.

Define
\begin{align}\label{correct}
\hat{P}(t)=\sup_{t-\tau(t)<s\leq t}\lVert\mu(s)\hat{Y}(s) \rVert_{\{\Lambda,\infty\}},
\end{align}
and with the same process in the above proof for Theorem 2, we have
\begin{align*}
		&\xi_{p_0}\frac{d\lVert \mu(t)\hat{Y}(t)\rVert_{\{\Lambda,\infty\}}}{dt}\bigg|_{t=t_0}\\
\leq&\bigg[\overline{T}_0\bigg(p_0,\frac{\dot{\mu}(t_0)}{\mu(t_0)}, \frac{\mu(t_0)}{\mu(t_0-\tau(t_0))}-1\bigg)-\kappa_{p_0}^{0}\Lambda_{p_0}[0]\bigg]\hat{P}(t_0)\\
&+\mu(t_0)(|\hat{I}_{p_0}^0|-\hat{\kappa}_{p_0}^0)\Lambda_p[0]<0.
		\end{align*}

Therefore, there exists a time point $T_1\ge T$ such that
\begin{align*}
		\sup_{t-\tau(t)\leq s \leq t}\lVert \hat{Y}(s) \rVert_{\{\Lambda, \infty\}} \leq 1, ~~~t\ge T_1
\end{align*}

\textbf{Phase II}: The norm of $\lVert \hat{Y}(t)\rVert_{\{\Lambda,\infty\}}$ defined in (\ref{normyt}) will be proved to evolve from 1 to 0 in finite time.
	
On the basis of (\ref{condition2}), we can find a sufficiently small constant $\theta>0$ such that, for $p=1,\cdots,n$, $\ell=0,1,\cdots,7$,
\begin{align}\label{useful4}
\sum_{q=1}^n \big|\tilde b_{pq}^T  M^{\ell}\big|  \overline{M}(q,\delta) \Lambda_q+|\hat{I}_p^{\ell}|-\hat{\kappa}_p^{\ell}+\Lambda_p[\ell]\theta<0.
\end{align}
	
	Define a Lyapunov function
	\begin{align}
		Y_2(t)=\lVert \hat{Y}(t) \rVert_{\{\Lambda,\infty\}}+\theta t,~~~t\ge T_1
	\end{align}
and the maximum-value function
	\begin{align}\label{SUPE2}
		P_2(t)=\sup_{t-\tau(t)\leq s\leq t}Y_2(s),~~~t\ge T_1.
	\end{align}
	
	If $Y_2(t)<P_2(t)$, then we can find a $\delta_2$-neighborhood such that $P_2(s)=P_2(t)$ for $s\in(t,t+\delta_2)$. Otherwise, $Y_2(t_0)=P_2(t_0)$, suppose $p_0=p_0(t_0)$ is an index which satisfies
	\begin{align}\label{p0-phase2}
		|\xi_{p_0}^{-1}\hat{w}_{p_0}^0(t_0)|+\theta t_0=P_2(t_0)
	\end{align}

Then simple calculations show that
	\begin{align*}
&\xi_{p_0}\frac{dP_2(t)}{dt}\bigg|_{t=t_0}
	=\frac{d}{dt}\big(|\hat{w}_{p_0}^0(t_0)|+\xi_{p_0}\theta t \big) \bigg  |_{t=t_0}\\
%	%%
%	\leq&\sum\limits_{q=1}^n|\tilde b_{p_0q}^T M^0|\overline{M}(q,\delta)|\tilde{e}_{q}(t_0-\tau_{p_0q}(t_0))|+|\hat{I}_{p_0}^0|-\hat{\kappa}_{p_0}^0+\xi_{p_0}\theta\\
%	%
%	%
	\leq& \sum\limits_{q=1}^n|\tilde b_{p_0q}^T M^0|\overline{M}(q,\delta)\Lambda_q+|\hat{I}_{p_0}^0|-\hat{\kappa}_{p_0}^0+\xi_{p_0}\theta<0
	\end{align*}

	Therefore, we can conclude that $P_2(t)$ is non-increasing for $t\ge T_1$, which leads to
	\begin{align*}
	Y_2(t)\leq P_2(t)\leq P_2(T_1)\le 1+\theta T_1
	\end{align*}
that is to say, when $T_2\ge T_1+1/\theta$, $\lVert \hat{Y}(t) \rVert_{\{\Lambda,\infty\}}=0$, i.e., OVNN would converge to the point $\hat{Z}^*$ in finite time.

\section*{Appendix D: Proof of Theorem \ref{theorem4}}
\textbf{Phase I}: The norm $\lVert \hat{Y}(t)\rVert_{\{\Lambda,\infty\}}$ defined in (\ref{normhyt}) is proved to decrease from initial value to 1 in finite time.
		
		Define a new Lyapunov function
		\begin{align}
			Y_3(t)=\mu(t)\lVert \hat{Y}(t) \rVert_{\{\Lambda,\infty\}}+\frac{1}{2c_2}(\kappa(t)-\kappa^*)^2
		\end{align}
		and denote
		\begin{align}\label{SUPE3}
			P_3(t)=\sup_{t-\tau(t)\leq s\leq t}Y_3(s).
		\end{align}

		If $Y_3(t)<P_3(t)$, then we can find a $\delta_3$, such that $P_3(s)=P_3(t)$ for $s\in(t,t+\delta_3)$. Otherwise,  there exists a $t_0\ge T$ such that $Y_3(t_0)=P_3(t_0)$, suppose index $p_0=p_0(t_0)$ satisfies
		\begin{align}
		\mu(t_0)|\xi_{p_0}^{-1}\hat{w}_{p_0}^0(t_0)|+\frac{1}{2c_2}(\kappa(t_0)-\kappa^*)^2=P_3(t_0)
		\end{align}

We will discuss firstly the limit of $\kappa(t)$. On the one hand, according to (\ref{rules of kappa}), it's obvious that $\kappa(t)$ is positive and monotonically increasing. Then based on the discussion in Theorem \ref{theorem3}, finite-time stability can be obtained if after some moment,
\begin{align}\label{fjx}
\kappa(t)\ge \kappa^*>\max_{p,\ell}\overline{T}_{\ell}(p, \alpha, \beta)\Lambda_p[\ell]^{-1}.
\end{align}
On the other hand, all the values of $\kappa(t)$ will be less than $\kappa^*$, function $(\kappa(t)-\kappa^*)^2$ is decreasing, thus $\mu(t)\lVert \hat{Y}(t) \rVert_{\{\Lambda,\infty\}}$ is the maximum value in $[t-\tau(t),t]$. Then,
		\begin{align*}
			&\xi_{p_0}\frac{dP_3(t)}{dt}=\frac{d}{dt} \big( \mu(t)|\hat{w}_{p_0}^0(t)|+\frac{\xi_{p_0}}{2c_2}(\kappa(t)-\kappa^*)^2 \big)\bigg|_{t=t_0}\\
			=&\dot{\mu}(t_0)|\hat{w}_{p_0}^0(t_0)|+\mu(t_0)\mathrm{sign}\big(\hat{w}_{p_0}^0(t_0)\big)\dot{\hat{w}}_{p_0}^0(t_0)\\
			&+(\kappa(t_0)-\kappa^*)\mu(t)|\hat{w}_{p_0}^0(t_0)|\\
			\leq& \bigg[\overline{T}_0\bigg(p_0,\frac{\dot{\mu}(t_0)}{\mu(t_0)}, \frac{\mu(t_0)}{\mu(t_0-\tau(t_0))}-1\bigg)-\kappa^*\Lambda_{p_0}[0]\bigg]\hat{P}(t_0)\\
&+\mu(t_0)(|\hat{I}_{p_0}^0|-\hat{\kappa}^*)\Lambda_{p_0}[0]<0.
		\end{align*}
		
Similar to the argument in Theorem \ref{theorem3}, there exists a time point $T_3\ge T$ such that
\begin{align*}
		\sup_{t-\tau(t)\leq s \leq t}\lVert \hat{Y}(s) \rVert_{\{\Lambda, \infty\}} \leq 1, ~~~t\ge T_3
\end{align*}

\textbf{Phase II}: According to (\ref{condition2}), we can find a sufficiently large constant $\hat{\kappa}^*>0$ and a small constant $\theta^*>0$, such that
\begin{align*}
		\sum_{q=1}^n \big|\tilde b_{pq}^T  M^{\ell}\big|  \overline{M}(q,\delta)\Lambda_q+|I_p^{\ell}|-\hat{\kappa}^*+\Lambda_p[\ell]\theta^*<0
		\end{align*}

Define
		\begin{align}
		Y_4(t)=&\lVert \hat{Y}(t) \rVert_{\{\Lambda,\infty\}}+\frac{1}{2c_1\min \Lambda}(\hat{\kappa}(t)-\hat{\kappa}^*)^2 \nonumber\\
		&+\frac{1}{2c_3}(\kappa(t)-\kappa^*)^2+\theta^*t
		\end{align}
and denote
		\begin{align}\label{SUPE4}
		P_4(t)=\sup_{t-\tau(t)\leq s\leq t}Y_4(s).
		\end{align}
		
The same argument can be applied on $\hat{\kappa}^*$ as that in Phase I, here we omit it. If there exists a time point $t_0\ge T_3$ such that $Y_4(t_0)=P_4(t_0)$, suppose index $p_0=p_0(t_0)$ satisfies
		\begin{align*}
		|\xi_{p_0}^{-1}\hat{w}_{p_0}^0(t_0)|&+\frac{1}{2c_1\min \Lambda}(\hat{\kappa}(t_0)-\hat{\kappa}^*)^2 \nonumber\\
			&+\frac{1}{2c_3}(\kappa(t_0)-\kappa^*)^2+\theta^*t_0=P_4(t_0)
		\end{align*}

Then, simple calculations as that in Theorem \ref{theorem3} show that
	\begin{align*}
		&\xi_{p_0}\frac{dP_4(t)}{dt}\bigg|_{t=t_0} \\
		=& \frac{d}{dt} \big( |\hat{w}_{p_0}^0(t)|+\frac{\xi_{p_0}}{2c_1\min \Lambda}(\hat{\kappa}(t_0)-\hat{\kappa}^*)^2\\
&+\frac{\xi_{p_0}}{2c_3}(\kappa(t_0)-\kappa^*)^2+\xi_{p_0}\theta^* t \big)\bigg|_{t=t_0}\\
		\le&|\dot{\hat{w}}_{p_0}^0(t_0)|+(\hat{\kappa}(t_0)-\hat{\kappa}^*)+(\kappa(t_0)-\kappa^*)|\hat{w}_{p_0}^0(t_0)|+\xi_{p_0}\theta^*\\
		\leq&\sum_{q=1}^n \big|\tilde b_{p_0q}^T  M^{\ell}\big|  \overline{M}(q,\delta)\Lambda_q+|\hat{I}_{p_0}^0|-\hat{\kappa}^*+\xi_{p_0}\theta^*<0
	\end{align*}
	
Therefore, $Y_4(t)\leq P_4(t)\leq P_4(T_3)\leq 1+\theta^* T_3+\frac{(\hat{\kappa}^*)^2}{2c_1\min \Lambda}+\frac{(\kappa^*)^2}{2c_3}$.
That is to say, when $T_4 \ge T_3+(1+\frac{(\hat{\kappa}^*)^2}{2c_1\min \Lambda}+\frac{(\kappa^*)^2}{2c_3})/\theta^*$, the norm $\lVert \hat{Y}(t) \rVert_{\{\Lambda,\infty\}}=0$, i.e., OVNN would converge to the given point $\hat{Z}^*$ in finite time. Moreover, adaptive control strengthes $\kappa(t)$ and $\hat{\kappa}(t)$ will also stop its increasing finally.

\end{document}